\documentclass[11pt]{article}
\usepackage{latexsym,amsfonts,euscript,amssymb}

\def\pf{\noindent{\bf Proof\quad }}
\def\pfend{\hfill{$\Box$}\vskip 0.2cm}

\newtheorem{lem}{\bf Lemma}[section]
\newtheorem{nt}[lem]{\bf Notation}
\newtheorem{defi}[lem]{\bf Definition}
\newtheorem{pro}[lem]{\bf Problem}
\newtheorem{thm}[lem]{\bf Theorem}
\newtheorem{rem}[lem]{\bf Remark}
\newtheorem{exa}[lem]{\bf Example}

\newtheorem{prop}[lem]{\bf Proposition}
\newtheorem{con}[lem]{\bf Conjecture}
{\left\lbrace\matrix}{endmatix\right\rbrace}

 {\title{ The properties of bordered matrix
 \linebreak of symmetric block design
\thanks{Project supported by the National Natural Science Foundation of China(Grant No. 11571121).}
}

\author{Mingchun Xu \thanks{E-mail: xumch@scnu.edu.cn}
\\  \small School of Mathematics,
  South-China Normal University, \\ Guangzhou, 510631, China
    }
\date{}

\begin{document}
\maketitle

\begin{abstract}   Let $X=\{p_1,p_2, \cdots, p_v\}$ be a $v$-set(a set of $v$  elements), called
points, and let $\mathcal{B}=\{ B_1, B_2,\cdots, B_v \}$ be a finite collection of subsets of
$X$, called blocks. The pair $(X, \mathcal{B})$ is called a  symmetric $(v, k,
\lambda)$  design if the following conditions hold:

(i)\,\, Each $B_i $ is a $k$-subset of  $X$.

(ii)\,\,\ Each $B_i \cap B_j$ is a $\lambda$-subset of $X$ for $ 1\leq i\neq j \leq v$.

(iii)\,\, The integers $v, k, \lambda$ satisfy  $0< \lambda <k <v-1$.
\begin{pro}
One of the major unsolved problems in combinatorics  is the determination of the precise range of values of $v, k ,$ and $\lambda$
for which a symmetric $(v,k,\lambda)$ design exists.
\end{pro}
The symmetric $(n^2+n+1, n+1,1)$ block design is a projective plane of order $n$.
Projective planes of order $n$ exist for all prime powers $n$ 
but for no
other $n$ is a construction known. Thus Conjecture 0.2 is the famous long-standing conjecture of finite  projective planes.
\begin{con} (O.Veblen and W.H.Bussey, 1906; R.H.Bruck and H.J. Ryser, 1949)\,\,  If a finite  projective plane of order $n$ exists, then $n$ is a power of some prime $p$.
\end{con}

It was proved in 1989 by a computer search that there does not exist any projective plane of order 10 by  Lam, C.W.H., Thiel, L. and Swiercz, S. Whether there exists any projective plane of order 12 is still open.

The author introduces  the bordered matrix of a  $(v,k,\lambda)$ symmetric  design, which  preserves  some row inner product property, and gives  some new necessary conditions for
the existence of the symmetric $(v,k,\lambda)$ design. Thus Theorem 0.3 generalizes  Schutzenberger's theorem and the Bruck-Ryser-Chowla
theorem on the existence of symmetric  block designs. This bordered matrix has been a breakthrough idea  since 1950.
\begin{thm}    Let  $C$ be  a $w$ by $w+d$ nonsquare rational  matrix without any column of $k\cdot 1^t_{w+d}$, where $1^t_{w+d}$ is the $w+d$-dimensional all 1 column vector and $k$ is a rational number, $\alpha, \beta$ be positive  integers
 such that  matrix $\alpha I_{w}+\beta J_{w}$ is the  positive definite matrix with  plus $d$  congruent factorization
$$ C\, C^t=\alpha I_{w}+\beta J_{w}.$$
Then the following cases hold.

\textbf{Case 1} If  $w\equiv 0\, (mod\, 4) $, $d=1$, then  $\beta $ is a perfect square.

\textbf{Case 2}  If $w\equiv 0\, (mod\, 4) $, $d=2$, then  $\beta $ is a sum of two squares.

  \textbf{Case 3}
If  $w\equiv 2\, (mod\, 4) $, $d=1$ and  $\alpha=a^2+b^2$, where $a, b$ are integers,  then  $\beta $ is a perfect square.

\textbf{Case 4}
If  $w\equiv 2\, (mod\, 4) $, $d=2$ and  $\alpha=a^2+b^2$, where $a, b$ are integers,  then  $\beta $ is a sum of two squares.

 \textbf{Case 5}
If  $w\equiv 1\, (mod\, 4) $, $d=1$ and  $\alpha=a^2+b^2$, where $a, b$ are integers,  then   $\alpha^*=\beta^* $, where $m^*$ denotes the square-free part of the integer $m$.

\textbf{Case 6}
If  $w\equiv 1\, (mod\, 4) $, $d=2$ and  $\alpha=a^2+b^2$, where $a, b$ are integers,  then     the equation
  $$ \alpha z^2= -x^2+ \beta y^2$$
must  have a solution in integers,
$x, y, z$,  not all zero.

 \textbf{Case 7}
If  $w\equiv 3\, (mod\, 4) $, $d=1$ ,  then   $\alpha^*=\beta^* $, where $m^*$ denotes the square-free part of the integer $m$.

\textbf{Case 8}
If  $w\equiv 3\, (mod\, 4) $, $d=2$ ,  then    the equation
  $$ \alpha z^2=- x^2+ \beta y^2$$
must  have a solution in integers,
$x, y, z$,  not all zero. \end{thm}
As an  application of the above theorem  the following theorems are obtained.
\begin{thm}    Conjecture 0.2 holds if  finite projective plane of order $n\leq 33$.\end{thm}
\begin{thm}   If $(v,k,\lambda)$ are parameters $(49,16,5),(154,18,2)$ and  $(115,19,3)$, then each symmetric $(v,k,\lambda)$ design does not exist.
\end{thm}

For Problem 0.1 or Conjecture 0.2 Lam's algorithm is an exponential time algorithm. But the proof of  main Theorem 0.3 is just the Ryser-Chowla
 elimination procedure. Thus author's algorithm is a polynomial time algorithm. It fully reflects  his algorithm high efficiency.

\end{abstract}


\textbf{2010 MR Subject Classification}\,\, 05B05, \, 11D09

\section{Introduction}
An incidence structure consists simply of a set $X$ of points and a set $\mathcal{B}$ of blocks, with a relation of incidence between points and blocks.
Symmetric  block designs have an enormous literature and discussions of their basic properties are readily available in \cite{h1,l5,ry}.

\begin{defi} Let $v, k$ and $\lambda$ be integers. 
  Let $X=\{p_1,p_2, \cdots, p_v\}$ be a $v$-set(a set of $v$  elements), called
points, and let $\mathcal{B}=\{ B_1, B_2,\cdots, B_v \}$ be a finite collection of subsets of
$X$, called blocks. The pair $(X, \mathcal{B})$ is called a  symmetric $(v, k,
\lambda)$  design if the following conditions hold:

(i)\,\, Each $B_i $ is a $k$-subset of  $X$.

(ii)\,\,\ Each $B_i \cap B_j$ is a $\lambda$-subset of $X$ for $ 1\leq i\neq j \leq v$.

(iii)\,\, The integers $v, k, \lambda$ satisfy  $0< \lambda <k <v-1$.

The set $\{ v, k, \lambda\}$ is called the set of parameters of the  symmetric
 $(X, \mathcal{B})$ design. We also use the notation $\mathcal{D}=(X,
\mathcal{B})$.
\end{defi}

\begin{defi} Define a $v\times v  $ 0-1 matrix $$A=(a_{i\, j})_{1\leq i\leq v,
1\leq j \leq v}, $$ whose rows are indexed by the points
$p_1,p_2,\cdots, p_v $ and columns are indexed by the blocks
$B_1,B_2,\cdots ,B_v$, by
\[\begin{array}{ccccc}
a_{i\, j}= \{
\begin{array}{ccccc} 1,\,\, \mbox{ if} \,\, p_i\in B_j,
\\ 0,\,\,  otherwise.
\end{array}
\end{array}
\]
Then $A$ is called the incidence matrix of the symmetric   $(v, k, \lambda)$ design.
We set $n=k-\lambda$ and call $n$ the order of the symmetric $( v, k, \lambda)$ design.
\end{defi}

$A^t$ denotes the transpose of $A$.
$J_v$ and $I_v$ are
the $v\times v$ all $1's$ matrix and the identity matrix,
respectively.

Let $A$ be a $v\times v$ 0-1 matrix. Then $A$ is the
incidence matrix of a symmetric $(v, k,\lambda)$ design if and only if 
$$ A\, A^t= \lambda J_v+(k-\lambda)I_v . $$

 In this paper we introduce  the bordered  matrix $C$ of the  $(v,k,\lambda)$ symmetric design and prove some new necessary conditions for the existence of the symmetric design.

 Then the bordered matrix $C$ of the  $(v,k,\lambda)$ symmetric design   is  obtained from $A$ by adding many rows rational vectors  and many  columns  rational vectors such that
 $$ C\, C^t= (\lambda+l) J_{v+s}+(k-\lambda)I_{v+s}  $$
 for some positive integers $l, s$, where $C$ is a nonsquare rational matrix and is of full row rank.

  Let $\alpha, \beta$ be positive  integers. A matrix $\alpha I_{w}+\beta J_{w}$ is called the  positive definite matrix with  plus $d$  congruent factorization property if  there exists  a nonsquare rational  $w$ by $w+d$ matrix  $C$  such that
$$ \alpha I_{w}+\beta J_{w}=C\, C^t, $$
where $d$ is the difference between the number of columns and the number of rows of $C$.

 In this paper we consider more the  positive definite matrix $\alpha I_{w}+\beta J_{w}$  with  the above plus $d$ congruent factorization property.

 What happens for   the  positive definite matrix $\alpha I_{w}+\beta J_{w}$  with  plus $d$  congruent factorization property
if   $\alpha, \beta$ are  two positive  integers?

If $d=1$ or $2$ then we obtain the following main theorem 1. As an application of the main theorem 1 it is easy to determine  that there does not exist finite projective plane of order $n$ if $n$ is each of
 the first open values 10, 12, 15, 18, 20, 24, 26 and 28, for which the Bruck-Ryser-Chowla Theorem can not  be used. For large $n$ the new method is also valid. Also some symmetric designs are excluded by the new method.

 Structure of the paper: some elementary definitions and results are  summarized and  the main theorems are stated  in \S 2 below.
  Some theorems from number theory  are needed  for our work in combinatorial analysis in \S 3. The proof of main theorem 1 will be given in \S 4.
  An application of main theorem 1  will be given in \S 5 and \S 6. Some concluding remarks will be given in \S 7.

\section{Background and statement of the main results  }
Symmetric  block designs have an enormous literature and discussions of their basic properties are readily available in \cite{h1,l5,ry}.

\begin{nt} $\mathbf{Z}$ denotes the set of integers.

$\mathbf{Q}$ denotes the field of rational numbers.

$m^*$ denotes the square-free part of the integer $m$.

Let $A$ be a matrix. $A^t$ denotes the transpose of $A$.

 $J_v$ and $I_v$ are
the $v\times v$ all $1's$ matrix and the identity matrix,
respectively.

$1_v$ is the $v$-dimensional all 1 row vector.
\end{nt}
\begin{prop} Let $A$ be a $v\times v$ 0-1 matrix. Then $A$ is the
incidence matrix of a symmetric $(v, k,\lambda)$ design if and only if 
\begin{equation}A\, A^t=A^t\, A= \lambda J_v+(k-\lambda)I_v .
\end{equation}
\end{prop}
\begin{prop}  In a symmetric  $(v,k, \lambda)$ design, the integers $v, k$, and $\lambda$ of the design must satisfy the following relations

(1) \,\   $\lambda(v-1)=k(k-1)$,

(2) \,\   $k^2-v\lambda =k-\lambda$, and

(3) \,\   $(v-k)\lambda =(k-1)(k-\lambda)$.
\end{prop}

\begin{thm} (Schutzenberger)\,\,  Suppose  there exists  a symmetric  $(v,k,\lambda)$ design with an incidence matrix $A$.
\,\ If $v$ is even, then $k-\lambda $ must be a perfect square.

\end{thm}
The Bruck-Ryser-Chowla Theorem gives a necessary condition for the existence of a symmetric design.

\begin{thm} (Bruck-Ryser-Chowla)\,\,  Suppose  there exists  a symmetric $(v,k,\lambda)$ design.
 If $v$ is odd, then  the equation
  $$ x^2= (k-\lambda)y^2+(-1)^{(v-1)/2}\lambda z^2$$
 must have a solution in integers,
$x, y, z$,  not all zero.
\end{thm}
\begin{rem}
 If $v$ is odd, then  the equation
  $$ x^2= (k-\lambda)y^2+(-1)^{(v-1)/2}\lambda z^2$$
 must have a solution in integers,
$x, y, z$,  not all zero. We also say that it has a nontrivial integral solution.
\end{rem}
\begin{pro}
One of the major unsolved problems in combinatories  is the determination of the precise range of values of $v, k ,$ and $\lambda$
for which a symmetric $(v,k,\lambda)$ design exists.
\end{pro}
\begin{thm} (Bruck-Ryser-Chowla)\,\,  \ Let  $C$ be  a square rational $w\times w$ matrix, $\alpha, \beta$ be positive integers such that
\begin{equation} C\, C^t=\alpha I_{w}+\beta J_{w}.\end{equation}
If $w$ is odd, then  the equation
  $$ z^2= \alpha x^2+(-1)^{(w-1)/2}\beta y^2$$
 must have a solution in integers,
$x, y, z$,  not all zero.
\end{thm}
\begin{defi} Let $n$ be a positive integer.   A finite projective plane of order $n$ is a symmetric  $(n^2+n, n+1, 1)$ design.  A block in a finite projective plane  is called a line.
\end{defi}
The theorem of Desargues is universally valid in a projective plane if and only if the plane can be constructed from a three-dimensional vector space over a field.  These planes are called Desarguesian planes, named after Girard Desargues.
The projective planes that can not be constructed in this manner are called non-Desarguesian planes, and the Moulton plane  is an example of one. The $PG(2,K)$ notation is reserved for the Desarguesian planes, where $K$ is some field.
\begin{thm}(O.Veblen and W.H.Bussey, 1906, see\cite{ve})
Let $q$ be a prime power. Then $PG(2, \mathbb{F}_q)$ is a finite projective plane of order of $q$.
\end{thm}
From Theorem 2.5 we deduce
\begin{thm} (Bruck-Ryser)\,\,  Let $n$ be a positive integer and $n\equiv \, 1, 2\, (mod\, 4)$ and let the squarefree part of $n$ contain at least one
prime factor $p\equiv \, 3\, (mod\, 4)$. Then there does not exist a finite  projective plane of order $n$.
\end{thm}
As an application, consider projective planes. Here $\lambda=1$ and $v=n^2+n+1$ is odd. If $n\equiv 0$ or $3\, (mod\,\ 4)$, the
Bruck-Ryser-Chowla equation always has the solution $(0,1,1)$ and thus the theorem excludes no values of $n$. However, if $n\equiv 1$ or $2\,\ (mod\,\  4)$,
 the equation becomes $nx^2=y^2+z^2$, which has a nontrivial integral solution if and only if $n$ is the sum of
 two squares of integers. Projective planes of order 6, 14, 21, 22, 30 or 33 therefore cannot exist.
\begin{con} (O.Veblen and W.H.Bussey, 1906; R.H.Bruck and H.J. Ryser, 1949)\,\,  If a finite  projective plane of order $n$ exists, then $n$ is a power of some prime $p$.
\end{con}

Despite much research no one has uncovered any further necessary conditions for the existence of a symmetric $(v,k,\lambda)$ design apart from
the equation $(v-1)\lambda=k(k-1)$, Schutzenberger's Theorem and the Bruck-Ryser-Chowla Theorem. For no $(v,k,\lambda)$ satisfying these requirements has it
been shown that a symmetric $(v,k,\lambda)$ design does not exist.

It is possible that these conditions are sufficient.  As a matter of fact, this is true the seventeen admissible $(v,k,\lambda)$ with $v\leq 48$(\cite{l5}, Notes to Chapter 2);
the first open case as of early 1982 is $(49, 16, 5)$.

Projective planes of order $n$ exist for all prime powers $n$ (aside from $PG(2,n)$ a host of other constructions are known ) but for no
other $n$ is a construction known. The first open values are $n=10, 12, 15, 18, 20, 24, 26$ and  28.  It was proved by a computer search that there does not exist any projective plane of order 10, cf. Lam, C.W.H., Thiel, L. and Swiercz, S. \cite{l4}. Whether there exists any projective plane of order 12 is still open.

Now we introduce  the bordered matrix of the  $(v,k,\lambda)$ symmetric design in Definition 2.14 and prove some new necessary conditions for the existence of the symmetric design.

\begin{rem}  The condition (1) in Proposition 2.2  for an incidence matrix $A$ of the symmetric  $(v,k,\lambda)$ design is equivalent to the following two conditions

 (i)\,\  the inner product of any two distinct rows of $A$ is equal to $\lambda$ ;

 (ii) \,\  and  the inner product of any rows with themselves of $A$ is equal to $k$.
  \end{rem}
  \begin{defi} Let $A$ be  an incidence matrix of the symmetric $(v,k,\lambda)$ design.
 Then the bordered matrix $C$ of $A$ for some positive integers $l$ and $d$ is  obtained from $A$ by adding many rows rational rectors
  and many columns  rational vectors such that

  (i)\,\, the inner product of any two distinct rows of $C$ is equal to $\lambda+l$ ;

  (ii)\,\, and
  the inner product of any rows with themselves of $C$ is equal to $k+l$;
\linebreak  where $C$ is  a $w$ by $w+d$ nonsquare rational matrix and is of full row rank.
 \end{defi}

\begin{rem} It is difficult to construct a square  bordered matrix of $A$. The author does this by computer computation in Maple.  But it is easy
 to construct a nonsquare   bordered  matrix of $A$. If it exists   then  it is not unique for some positive integers $l$ and $d$. The author also does this by computer computation in Maple.
 \end{rem}
 Theorem 2.8 gives a  necessary condition for the existence of positive definite matrix $\alpha I_{w}+\beta J_{w}$, which is congruent to identity matrix over rational field for the positive integers $\alpha, \beta$.
  \begin{defi}  Let $\alpha, \beta$ be positive  integers. A matrix $\alpha I_{w}+\beta J_{w}$ is called the  positive definite matrix with  plus $d$  congruent factorization property if  there exists  a nonsquare rational  $w$ by $w+d$ matrix  $C$  such that
\begin{equation} \alpha I_{w}+\beta J_{w}=C\, C^t.\end{equation}
\end{defi}
\begin{rem} The matrix equation (3) implies $C$ is  always of rank $w$, i.e., of full row rank if   $\alpha, \beta$ are  two positive  integers and $C$ is a
 $w$ by $w+d$ matrix over rational field $\mathbf{Q}$.
 \end{rem}
 In this paper we consider more the  positive definite matrix $\alpha I_{w}+\beta J_{w}$  with  plus $d$ congruent factorization property.
\begin{pro} What happens for   the  positive definite matrix $\alpha I_{w}+\beta J_{w}$  with  plus $d$  congruent factorization property
if   $\alpha, \beta$ are  two positive  integers?

\end{pro}
The matrix equation (3) is of fundamental importance. But it is difficult
to deal with this matrix equation in its full generality. Nevertheless,  if $d=1$ or $2$ then we obtain the following main theorem.
Thus  Main Theorem 1 generalizes  Schutzenberger's theorem and the Bruck-Ryser-Chowla
theorem on the existence of symmetric  block designs.

We are now prepared to state  our main conclusions.

\textbf{ Main Theorem 1}   Let  $C$ be  a $w$ by $w+d$ nonsquare rational  matrix without any column of $k\cdot 1^t_{w+d}$, where $1^t_{w+d}$ is the $w+d$-dimensional all 1 column vector and $k$ is a rational number, $\alpha, \beta$ be positive  integers
 such that  matrix $\alpha I_{w}+\beta J_{w}$ is the  positive definite matrix with  plus $d$  congruent factorization
$$ C\, C^t=\alpha I_{w}+\beta J_{w}.$$
Then the following cases hold.

\textbf{Case 1} If  $w\equiv 0\, (mod\, 4) $, $d=1$, then  $\beta $ is a perfect square.

\textbf{Case 2}  If $w\equiv 0\, (mod\, 4) $, $d=2$, then  $\beta $ is a sum of two squares.

  \textbf{Case 3}
If  $w\equiv 2\, (mod\, 4) $, $d=1$ and  $\alpha=a^2+b^2$, where $a, b$ are integers,  then  $\beta $ is a perfect square.

\textbf{Case 4}
If  $w\equiv 2\, (mod\, 4) $, $d=2$ and  $\alpha=a^2+b^2$, where $a, b$ are integers,  then  $\beta $ is a sum of two squares.

 \textbf{Case 5}
If  $w\equiv 1\, (mod\, 4) $, $d=1$ and  $\alpha=a^2+b^2$, where $a, b$ are integers,  then   $\alpha^*=\beta^* $, where $m^*$ denotes the square-free part of the integer $m$.

\textbf{Case 6}
If  $w\equiv 1\, (mod\, 4) $, $d=2$ and  $\alpha=a^2+b^2$, where $a, b$ are integers,  then     the equation
  $$ \alpha z^2= -x^2+ \beta y^2$$
must  have a solution in integers,
$x, y, z$,  not all zero.

 \textbf{Case 7}
If  $w\equiv 3\, (mod\, 4) $, $d=1$ ,  then   $\alpha^*=\beta^* $, where $m^*$ denotes the square-free part of the integer $m$.

\textbf{Case 8}
If  $w\equiv 3\, (mod\, 4) $, $d=2$ ,  then    the equation
  $$ \alpha z^2=- x^2+ \beta y^2$$
must  have a solution in integers,
$x, y, z$,  not all zero.

As an application of the main theorems it is easy to determine  that there does not exist finite projective plane of order $n$ if $n$ is each of
 the first open values 10, 12, 15, 18, 20, 24, 26 and 28, for which the Bruck-Ryser-Chowla Theorem can not  be used. For large $n$ the new method is also valid. Also some symmetric designs are excluded by the new method.

\textbf{ Main Theorem 2}   Conjecture 2.12 holds if  finite projective plane of order $n\leq 33$.

\textbf{ Main Theorem 3}  If $(v,k,\lambda)$ are paremeters $(49,16,5),(154,18,2)$ and $(115,19,3)$, then the  symmetric $(v,k,\lambda)$ designs do not exist.

\section{Some theorems from number theory}
In this section we shall state some theorems from number theory that are needed  for our work in combinatorial analysis. No proofs will be given, but references will
be given to books where the proofs may be found.

\begin{lem}\,\  (Lagrange,  Sum of Four Squares Theorem, \cite{se}) Every positive integer is the sum of four integral squares.
 \end{lem}
 \begin{lem} ( Sum of Two Squares Theorem, \cite{si}, Theorem 27.1) Let $m$ be a positive integer. Factor $m$ as $$m=p_1p_2\cdots p_rM^2$$ with distinct prime factors
 $p_1,p_2,\cdots,p_r $. Then $m$ can be written as a sum of two integral squares exactly when  each   $p_i$ is either 2 or is congruent to 1 $modulo$ 4.
 \end{lem}
\begin{lem}  ( Sum of Three Squares Theorem, \cite{se}) Positive integer $n$ is the sum of three  integral squares if  $n^* \equiv 1, 2, 3, 5, $ or $ 6\,\  (mod\,\ 8) $.
 \end{lem}
We shall also use the following elementary identities which can be verified by direct multiplication.
\begin{lem}  (  Two Squares Identity, \cite{si}, Chapter 26)\, \,  $(b_1^2+b_2^2)(x_1^2+x_2^2)=y_1^2+y_2^2$, where
$$y_1=b_1x_1-b_2x_2 ,$$
$$y_2=b_2x_1+b_1x_2.$$
 \end{lem}
\begin{lem}  ( Four Squares Identity, \cite{l5}, \S 2.1)\, \,  $(b_1^2+b_2^2+b_3^2+b_4^2)(x_1^2+x_2^2+x_3^2+x_4^2)=y_1^2+y_2^2+y_3^2+y_4^2$, where
$$y_1=b_1x_1-b_2x_2-b_3x_3-b_4x_4,$$
$$y_2=b_2x_1+b_1x_2-b_4x_3+b_3x_4,$$
$$y_3=b_3x_1+b_4x_2+b_1x_3-b_2x_4,$$
$$y_4=b_4x_1-b_3x_2+b_2x_3+b_1x_4.$$
 \end{lem}

Let $m$ and $b$ be nonzero two integers and $(b, m)=1$. The integers $b$ are divided into two classes called \textit{quadratic residues} and
\textit{quadratic nonresidues} according as  $x^2\equiv b(mod\, m)$ does or does not have a solution $x(mod\, m)$.

Let $p$ be an odd prime. The integers $b$ with $p \nmid b $ are divided into two classes called quadratic residues and  quadratic nonresidues
according as $x^2\equiv b(mod\, p)$ does or does not have a solution $x(mod\, p)$. This property is expressed in term of the $Legendre\, symbol\, (\frac{b}{p})$ by
the rules
$$ (\frac{b}{p})=+1\,\ \mbox{if\, b\, a\, quadratic\,  residue \, modulo\, p, } $$
$$ (\frac{b}{p})=-1\,\ \mbox{if\, b\, a\, quadratic\,  nonresidue \, modulo\, p. } $$

\begin{lem}\,\, (Legendre, \cite{l5}, \S 2.1) Let $a, b, c$ be all positive, coprime to each other  and square-free integers.  The equation
\begin{equation} ax^2+by^2=cz^2\end{equation}
has solutions in integers $x, y, z$ not all zero if and only if $bc, ac$ and  $-ab$ are quadratic residues $mod(a)$,  $mod(b)$ and $mod(c)$ respectively.

 \end{lem}
\begin{lem} (Legendre, \cite{l5}, \S 2.1) Consider the equation \begin{equation} Ax^2+By^2+Cz^2=0 ,\end{equation}
 and assume initially that $A,B$, and $C$ are square-free integers, pairwise relatively prime. Necessary conditions for the existence of a
 nontrivial integral solution are that, for all odd primes $p$,

 (1)\,\ If $p|\, A$, then the Legendre symbol $(\frac{-BC}{p})=1$,

 (2)\,\ If $p|\, B$, then the Legendre symbol $(\frac{-AC}{p})=1$,

 (3)\,\ If $p|\, C$, then the Legendre symbol $(\frac{-AB}{p})=1$,

 and, of course,

 (4)\,\ $A,B$, and $C$ do not all have the same sign.

 It is a classical theorem, due to  Legendre that these simple necessary conditions are sufficient.\end{lem}

 \begin{rem} (\cite{l5}, \S 2.1) \, \, If  $A,B$, and $C$ do not satisfy our assumptions above we may slightly modify the equation (5). Henceforth, let $m^{*}$  denote the square-free
 part of the integer $m$. Then (5) has a
 nontrivial integral solution  if and only if \begin{equation} A^*x^2+B^*y^2+C^*z^2=0 ,\end{equation}
 has a
 nontrivial integral solution. \end{rem}
 \begin{rem} (\cite{l5}, \S 2.1 )\,\,
 If $p$ divides all three coefficients; we may divide it out and if $p$ divides only $A$ and $B$, then (5) has a nontrivial integral solution if
 only if $$ \frac{A}{p}x^2+\frac{B}{p}y^2+(pC)z^2=0 $$ does. Hence (6) can always be transformed into an equation to which Legendre's result applies.
\end{rem}

\section{Proof of   main Theorem 1}

 We are going to complete the proof of  main theorem 1 in this section and to give some examples of the existence of
 bordered matrix of symmetric $(v,k,\lambda)$ designs. We will see that the proof of  main Theorem 1 is just like  the Ryser-Chowla
 elimination procedure  in \cite{ch}. Also it is just like the Gausssian elimination procedure for solving the homogeneous linear equations.  

\begin{lem} (\textbf{Case 1 of Main Theorem 1} )  Let  $C$ be  a $w$ by $w+1$ nonsquare rational  matrix without any column of $k\cdot 1^t_{w+1}$, where $1^t_{w+1}$ is the $w+1$-dimensional all 1 column vector and $k$ is a rational number, $\alpha, \beta$ be two positive  integers. Suppose
  the matrix $\alpha I_{w}+\beta J_{w}$ is the  positive definite matrix with  plus $1$  congruent factorization
 property such that
\begin{equation} C\, C^t=\alpha I_{w}+\beta J_{w}.\end{equation}
If  $w\equiv 0\, (mod\, 4) $, then  $\beta $ is a perfect square.
\end{lem}
\pf
By the assumption   we have the identity
$$
 C \, \, C^t= \alpha I_{w}+  \beta  J_{w}$$ for the rational  matrix $C$.
 The idea of the proof is to interpret this as an identity in quadratic forms over the rational field.

  Suppose that $w\equiv 0\, (mod\, 4)$.
  If $x$ is the row vector $(x_1,x_2,\cdots, x_{w})$, then the
  identity for $C\, C^t$ gives $$xC\, C^t x^t= \alpha(x_1^2+x_2^2+\cdots+x_{w}^2)+\beta(x_1+x_2+\cdots+x_{w})^2. $$
Putting $f=x\, C$, $f=(f_1,f_2,\cdots, f_{w}, f_{w+1})$, $z=x_1+x_2+\cdots+x_{w}$, we have $f\, f^t= x C\, C^tx^t$ and  
\begin{equation} f_1^2+f_2^2+\cdots+f_{w}^2+f_{w+1}^2
=
\alpha (x_1^2+x_2^2+\cdots+x_{w}^2)+\beta z ^2;
\end{equation}
$$f_1=c_{1\, 1}x_1+c_{2\, 1}x_2\cdots+c_{w\, 1} x_{w},  $$
$$f_2=c_{1\, 2}x_1+c_{2\, 2}x_2\cdots+c_{w\, 2} x_{w},  $$
$$\cdots \cdots\cdots\cdots\cdots   $$
$$f_w=c_{1\, w}x_1+c_{2\, w}x_2\cdots+c_{w\, w} x_{w},  $$
$$f_{w+1}=c_{1\, {w+1}}x_1+c_{2\, {w+1}}x_2\cdots+c_{w\, {w+1}} x_{w};  $$
$$z=x_1+x_2+\cdots+x_{w} .$$
Thus the cone (8) of variables $f_1, f_2, \cdots, f_{w}, f_{w+1}, x_1,x_2, \cdots, x_{w}, z$ has some nontrivial rational points. We will get a nontrivial rational point for $y^2=\beta z^2$ such that $y\neq 0$ by the Ryser-Chowla elimination procedure for the above homogeneous equations.

Now define a linear mapping $\sigma$ from $\mathbf{Q}^{w}$ to $\mathbf{Q}^{w+1}$
$$ \sigma : x\mapsto x\ C.$$
The image space $\sigma(\mathbf{Q}^w)$ is a vector subspace of $\mathbf{Q}^{w+1}$. Let $\gamma_1, \gamma_2, \cdots, \gamma_w$ be the row vectors  of $C$. Thus the row space $R(C)$ is subspace of $\mathbf{Q}^{w+1}$ spanned  by $\gamma_1, \gamma_2, \cdots, \gamma_w$. So $\sigma(\mathbf{Q}^w)=R(C)$. By Remark 2.17, since  $C$ is of full row rank,  $dim_{\mathbf{Q}}(R(C))=w.$ So $\sigma$ is an one-one linear mapping  from $\mathbf{Q}^w$ to $R(C)$.

The equation (8) is an identity in $x_1, x_2, \cdots,x_{w} .$
Each of the f's is a rational combination of the x's, since $f=x\, C$. By Remark 2.17, since  $C$ is of full row rank, each of the x's is a rational combination of the f's for any  $f\in R(C)$. Thus the equation (8) is an identity in
 the variables $f_1, f_2, \cdots,f_{w}, f_{w+1}$ for any  $f\in R(C)$.

 We express the integer $\alpha$ as the sum of four squares by Lemma 3.1, and bracket the terms $x_1^2+\cdots +x_{w}^2$  in fours.
 Each product of sums of four squares is itself a sum of four squares, i.e. Lemma 3.5,  and  so (8) yields $$f_1^2+f_2^2+\cdots + f_{w}^2+f_{w+1}^2$$
\begin{equation}
 = y_1^2+y_2^2+\cdots+y_{w}^2+ \beta  z^2,
\end{equation}
where $z=x_1+x_2+\cdots+x_{w}$, and the $y's$ are related to the
$x's$ by an invertible linear transformation with rational
coefficients. Since the $y's$ are rational linear combinations of
the $x's$, it follows that the $y's $ (and $z$) are rational linear
combinations of the $f's$ for any  $f\in R(C)$.  Thus the equation (9) is an identity in
 the variables $f_1, f_2, \cdots,f_{w}, f_{w+1}$ for any  $f\in R(C)$.

Suppose that $ y_i=b_{i\, 1}f_1+\cdots +b_{i\, w}f_{w}+b_{i\, w+1}f_{w+1},\, 1 \leq i \leq w $. We can
define $f_1$ as a rational linear combination of $f_2,\cdots ,
f_{w+1}$, in such a way that $y_1^2=f_1^2$: if $b_{1\, 1}\neq 1$ we
set $f_1=\frac{1}{1-b_{1\, 1}}(b_{1\, 2}f_2+\cdots+b_{1\,
w+1}f_{w+1})$, while if  $b_{1\, 1}= 1$ we set
$f_1=\frac{1}{-1-b_{1\, 1}}(b_{1\, 2}f_2+\cdots+b_{1\,
w+1}f_{w+1})$. Now we know that  $y_2$ is  a rational linear
combination of the $f's$, and, using the relevant expression for
$f_1$ found above, we can express $y_2$  as a rational linear
combination of $f_2,\cdots,f_{w+1}$. As before, we fix $f_2$ as a
rational combination of $f_3,\cdots , f_{w+1}$ in such a way that
$y_2^2=f_2^2$. Continuing thus, we eventually obtain $y_1,\cdots
,y_{w}$ and $f_1,\cdots ,f_{w}$ as rational multiples  of
$f_{w+1}$, satisfying $f_i^2=y_i^2\, (1 \leq i \leq w)$.

We reduce the equations step by step in this way until a truncated triangle of equations is obtained, say
 $$f_1=d_{1\, 2}f_2+\cdots+d_{1\,
w+1}f_{w+1},  $$
$$f_2=d_{2\, 3}f_3+\cdots+d_{2\,
w+1}f_{w+1},  $$
$$\cdots \cdots\cdots\cdots\cdots   $$
$$f_w=d_{w\, \, w+1}f_{w+1};  $$
$$f_i^2=y_i^2\, , (1 \leq i \leq w) ; $$
where $d_{i\, j} \in \mathbf{ Q}$.

For any  $x\in \mathbf{Q}^w$ and $x\neq 0$,  by Remark 2.17, since  $C$ is of full row rank, $f=x\, C$ implies $f\neq 0$. Let the last one $f_{w+1}\neq 0$ with suitable renumberings, if necessary.
Choose any non-zero rational value for  $f_{w+1} $.
All the
$y's $, the remaining $f's$, and $z$, are determined as above, and
substituting these values in (9) we obtain
\begin{equation}f_{w+1}^2= \beta z^2.
\end{equation}
Multiplying by a suitable constant  we have  that $\beta $ is a perfect square. So the theorem is proved.  \pfend

\begin{exa}
 The projective plane of order 5 is the symmetric $(31,6, 1)$ design.
Let $A$ be its incidence matrix, which is  a 31 by 31 matrix.
Choose its  bordered matrix $C$ is  a 32 by 33 matrix
as the following matrix.
 \[\begin{array}{ccccc} C=\left[\begin{array}{ccccc} A_{1\,\, 1} \ A_{1\,\, 2}
\\ A_{2\,\, 1}\ A_{2\,\, 2}
\end{array}\right ],
\end{array}\]
$$A_{1\,\, 1}=A .$$ $A_{1\,\, 2}$ is a 31 by 2 matrix  and
$$ A_{1\,\, 2}
=(2*1_{31}^t,  2*1_{31}^t). $$
 $A_{2\,\, 1}$ is
 a 1 by 31 matrix and
  \[\begin{array}{ccccc} A_{2\,\, 1}=\left[\begin{array}{ccccc}    \frac{1}{12}\cdot1_{31}
\end{array}\right ].
\end{array}\]
$A_{2\,\, 2}$ is a 1 by 2 matrix and
 \[ \begin{array}{ccccc} A_{2\,\, 2}=\left[\begin{array}{ccccccccccccccc}   \frac{7}{12} \ \frac{11}{3}   \end{array}\right ] . \end{array}
\]
  It is easy to check that $C$ has the property of row inner products, i.e.,

  (i)\,\, the inner product of any two distinct rows of $C$ is equal to $9$ ;

  (ii)\,\, and
  the inner product of any rows with themselves of $C$ is equal to $14$.\linebreak
  It follows that
 $$ C\, C^t=5 I_{32}+9 J_{32}$$
  and $C$ is exactly  the  bordered matrix of the symmetric $(31,6, 1)$ design. You wish to verify this by hand or electronic computation.
 We have that 9  is  a  perfect square just as the assertion of the above theorem.
\end{exa}

\begin{lem} (\textbf{ Case 2 of Main Theorem 1} )  Let  $C$ be  a $w$ by $w+2$ nonsquare rational  matrix without any column of $k\cdot 1^t_{w+2}$, where $1^t_{w+2}$ is the $w+2$-dimensional all 1 column vector and $k$ is a rational number, $\alpha, \beta$ be positive  integers.
 Suppose the  matrix $\alpha I_{w}+\beta J_{w}$ is the  positive definite matrix with  plus $2$  congruent factorization
 property such that
\begin{equation} C\, C^t=\alpha I_{w}+\beta J_{w}.\end{equation}
If $w\equiv 0\, (mod\, 4) $, then  $\beta $ is a sum of two squares.
\end{lem}
\pf
By the assumption   we have the identity
$$
 C \, \, C^t= \alpha I_{w}+  \beta  J_{w}$$ for the rational  matrix $C$.
 The idea of the proof is to interpret this as an identity in quadratic forms over the rational field.

  Suppose that  $w$ is an even integer with $w\equiv 0\, (mod\, 4)$.
  If $x$ is the row vector $(x_1,x_2,\cdots, x_{w})$, then the
  identity for $C\, C^t$ gives $$xC\, C^t x^t= \alpha(x_1^2+x_2^2+\cdots+x_{w}^2)+\beta(x_1+x_2+\cdots+x_{w})^2. $$
Putting $f=x\, C$, $f=(f_1,f_2,\cdots, f_{w}, f_{w+1},  f_{w+2})$, $z=x_1+x_2+\cdots+x_{w}$, we have $f\, f^t= x C\, C^tx^t$ and \begin{equation} f_1^2+f_2^2+\cdots+f_{w}^2+f_{w+1}^2 +f_{w+2}^2
=
\alpha(x_1^2+x_2^2+\cdots+x_{w}^2)+\beta z^2;
\end{equation}
$$f_1=c_{1\, 1}x_1+c_{2\, 1}x_2\cdots+c_{w\, 1} x_{w},  $$
$$f_2=c_{1\, 2}x_1+c_{2\, 2}x_2\cdots+c_{w\, 2} x_{w},  $$
$$\cdots \cdots\cdots\cdots\cdots   $$
$$f_w=c_{1\, w}x_1+c_{2\, w}x_2\cdots+c_{w\, w} x_{w},  $$ $$f_{w+1}=c_{1\, {w+1}}x_1+c_{2\, {w+1}}x_2\cdots+c_{w\, {w+1}} x_{w},  $$
$$f_{w+2}=c_{1\, {w+2}}x_1+c_{2\, {w+2}}x_2\cdots+c_{w\, {w+2}} x_{w};  $$
$$z=x_1+x_2+\cdots+x_{w} .$$
Thus the cone (12) of variables $f_1, f_2, \cdots, f_{w}, f_{w+1},f_{w+2}, x_1,x_2, \cdots, x_{w}, z$ has some nontrivial rational points. We will get a nontrivial rational point for $f_{w+1}^2+f_{w+2}^2= \beta z^2$  such that $f_{w+2}^2 \neq 0$ by the Ryser-Chowla elimination procedure for the above homogeneous equations.

Now define a linear mapping $\sigma$ from $\mathbf{Q}^{w}$ to $\mathbf{Q}^{w+2}$
$$ \sigma : x\mapsto x\ C.$$
The image space $\sigma(\mathbf{Q}^w)$ is a vector subspace of $\mathbf{Q}^{w+2}$. Let $\gamma_1, \gamma_2, \cdots, \gamma_w$ be the row vectors  of $C$. Thus the row space $R(C)$ is subspace of $\mathbf{Q}^{w+2}$ spanned  by $\gamma_1, \gamma_2, \cdots, \gamma_w$. So $\sigma(\mathbf{Q}^w)=R(C)$. By Remark 2.17, since  $C$ is of full row rank,  $dim_{\mathbf{Q}}(R(C))=w.$ So $\sigma$ is an one-one linear mapping  from $\mathbf{Q}^w$ to $R(C)$.

The equation (12) is an identity in $x_1, x_2, \cdots,x_{w} .$
Each of the f's is a rational combination of the x's, since $f=x\, C$. By Remark 2.17, since  $C$ is of full row rank, each of the x's is a rational combination of the f's. Thus the equation (12) is an identity in
 the variables $f_1, f_2, \cdots,f_{w}, f_{w+1}, f_{w+2}$ for any  $f\in R(C)$.

 We express the integer $\alpha$ as the sum of four squares by Lemma 3.1, and bracket the terms $x_1^2+\cdots +x_{w}^2$  in fours.
 Each product of sums of four squares is itself a sum of four squares, i.e. Lemma 3.5, and  so (12) yields $$f_1^2+f_2^2+\cdots + f_{w}^2+f_{w+1}^2+f_{w+2}^2$$
\begin{equation}
 = y_1^2+y_2^2+\cdots+y_{w}^2+ \beta  z^2,
\end{equation}
where $z=x_1+x_2+\cdots+x_{w}$, and the $y's$ are related to the
$x's$ by an invertible linear transformation with rational
coefficients. Since the $y's$ are rational linear combinations of
the $x's$, it follows that the $y's $ (and $z$) are rational linear
combinations of the $f's$.  Thus the equation (13) is an identity in
 the variables $f_1, f_2, \cdots,f_{w}, f_{w+1}, f_{w+2}$ for any  $f\in R(C)$.

Suppose that   $ y_i=b_{i\, 1}f_1+\cdots +b_{i\, w}f_{w}+b_{i\, w+1}f_{w+1} +b_{i\, w+2}f_{w+2}, \, 1 \leq i \leq w $.  We can
define $f_1$ as a rational linear combination of $f_2,\cdots ,
f_{w+1}, f_{w+2}$, in such a way that $y_1^2=f_1^2$: if $b_{1\, 1}\neq 1$ we
set $f_1=\frac{1}{1-b_{1\, 1}}(b_{1\, 2}f_2+\cdots+b_{1\,
w+1}f_{w+1 +b_{1\, w+2}f_{w+2}})$, while if  $b_{1\, 1}= 1$ we set
$f_1=\frac{1}{-1-b_{1\, 1}}(b_{1\, 2}f_2+\cdots+b_{1\,
w+1}f_{w+1} +b_{1\, w+2}f_{w+2} )$. Now we know that  $y_2$ is  a rational linear
combination of the $f's$, and, using the relevant expression for
$f_1$ found above, we can express $y_2$  as a rational linear
combination of $f_2,\cdots,f_{w+1},f_{w+2}$. As before, we fix $f_2$ as a
rational combination of $f_3,\cdots , f_{w+1}, f_{w+2}$ in such a way that
$y_2^2=f_2^2$. Continuing thus, we eventually obtain $y_1,\cdots
,y_{w}$ and $f_1,\cdots ,f_{w}$ as rational combinations  of
$f_{w+1}, f_{w+2}$, satisfying $f_i^2=y_i^2\, (1 \leq i \leq w)$.

We reduce the equations step by step in this way until a truncated triangle of equations is obtained, say
 $$f_1=d_{1\, 2}f_2+\cdots+d_{1\,
w+1}f_{w+1}+d_{1\,
w+2}f_{w+2},  $$
$$f_2=d_{2\, 3}f_3+\cdots+d_{2\,
w+1}f_{w+1}+d_{2\,
w+2}f_{w+2},  $$
$$\cdots \cdots\cdots\cdots\cdots   $$
$$f_w=d_{w\, \, w+1}f_{w+1}+d_{w\,
w+2}f_{w+2} ;$$
$$f_i^2=y_i^2, \, (1 \leq i \leq w); $$
where $d_{i\, j} \in \mathbf{ Q}$.

For any  $x\in \mathbf{Q}^w$ and $x\neq 0$,  by Remark 2.17, since  $C$ is of full row rank, $f=x\, C$ implies $f\neq 0$. Let the last one $f_{w+2}\neq 0$ with suitable renumberings, if necessary.
Choose any non-zero rational value for $f_{w+2} $.
All the
$y's $, the remaining $f's$, and $z$, are determined as above, and
substituting these values in (13) we obtain
\begin{equation}f_{w+1}^2+f_{w+2}^2= \beta z^2.
\end{equation}
Multiplying by a suitable constant  we have  that $\beta $ is a sum of two squares. So the theorem is proved.  \pfend

\begin{exa}
 The projective plane of order 5 is the symmetric $(31,6, 1)$ design.
Let $A$ be its incidence matrix, which is  a 31 by 31 matrix.
Choose its  bordered matrix $C$ is  a 32 by 34 matrix
as the following matrix.
 \[\begin{array}{ccccc} C=\left[\begin{array}{ccccc} A_{1\,\, 1} \ A_{1\,\, 2}
\\ A_{2\,\, 1}\ A_{2\,\, 2}
\end{array}\right ],
\end{array}\]
$$A_{1\,\, 1}=A .$$ $A_{1\,\, 2}$ is a 31 by 3 matrix  and
$$ A_{1\,\, 2}
=(1*1_{31}^t,  0*1_{31}^t,0*1_{31}^t). $$
 $A_{2\,\, 1}$ is
 a 1 by 31 matrix and
  \[\begin{array}{ccccc} A_{2\,\, 1}=\left[\begin{array}{ccccc}    \frac{1}{3}\cdot1_{31}
\end{array}\right ].
\end{array}\]
$A_{2\,\, 2}$ is a 1 by 3 matrix and
 \[ \begin{array}{ccccc} A_{2\,\, 2}=\left[\begin{array}{ccccccccccccccc}  0 \  \frac{4}{3} \ \frac{4}{3}   \end{array}\right ] . \end{array}
\]
  It is easy to check that $C$ has the property of  row inner products, i.e.,

  (i)\,\, the inner product of any two distinct rows of $C$ is equal to $2$ ;

  (ii)\,\, and
  the inner product of any rows with themselves of $C$ is equal to $7$.\linebreak
  It follows that
 $$ C\, C^t=5 I_{32}+ 2 J_{32}$$
  and $C$ is exactly  the  bordered matrix of the symmetric $(31,6, 1)$ design.
 We have that 2  is  is a sum of two squares  just as the assertion of the above theorem.
\end{exa}

\begin{lem}  (\textbf{Case 3 of Main Theorem 1} ) Let  $C$ be  a $w$ by $w+1$ nonsquare rational  matrix without any column of $k\cdot 1^t_{w+1}$, where $1^t_{w+1}$ is the $w+1$-dimensional all 1 column vector and $k$ is a rational number,
$\alpha, \beta$ be positive  integers and
 $\alpha=a^2+b^2$, where $a, b$ are   integers. Suppose  the
 matrix $\alpha I_{w}+\beta J_{w}$ is the  positive definite matrix with  plus $1$  congruent factorization
 property such that
\begin{equation} C\, C^t=\alpha I_{w}+\beta J_{w}.\end{equation}
If  $w\equiv 2\, (mod\, 4) $, then  $\beta $ is a perfect square.
\end{lem}
\pf
By the assumption   we have the identity
$$
 C \, \, C^t= \alpha I_{w}+  \beta  J_{w}$$ for the rational  matrix $C$.
 The idea of the proof is to interpret this as an identity in quadratic forms over the rational field.

  Suppose that   $w\equiv 2\, (mod\, 4)$.
  If $x$ is the row vector $(x_1,x_2,\cdots, x_{w})$, then the
  identity for $C\, C^t$ gives $$xC\, C^t x^t= \alpha(x_1^2+x_2^2+\cdots+x_{w}^2)+\beta(x_1+x_2+\cdots+x_{w})^2. $$
Putting $f=x\, C$, $f=(f_1,f_2,\cdots, f_{w}, f_{w+1})$, $z=x_1+x_2+\cdots+x_{w}$, we have $f\, f^t= x C\, C^tx^t$ and  
\begin{equation}f_1^2+f_2^2+\cdots+f_{w}^2+f_{w+1}^2
=
\alpha(x_1^2+x_2^2+\cdots+x_{w}^2)+\beta z^2;
\end{equation}
$$f_1=c_{1\, 1}x_1+c_{2\, 1}x_2\cdots+c_{w\, 1} x_{w},  $$
$$f_2=c_{1\, 2}x_1+c_{2\, 2}x_2\cdots+c_{w\, 2} x_{w},  $$
$$\cdots \cdots\cdots\cdots\cdots   $$
$$f_w=c_{1\, w}x_1+c_{2\, w}x_2\cdots+c_{w\, w} x_{w},  $$
$$f_{w+1}=c_{1\, {w+1}}x_1+c_{2\, {w+1}}x_2\cdots+c_{w\, {w+1}} x_{w};  $$
$$z=x_1+x_2+\cdots+x_{w} .$$
Thus the cone (16) of variables $f_1, f_2, \cdots, f_{w}, f_{w+1}, x_1,x_2, \cdots, x_{w}, z$ has some nontrivial rational points. We will get a nontrivial rational point for $y^2=\beta z^2$ such that $y\neq 0$ by the Ryser-Chowla elimination procedure for the above homogeneous equations.

Now define a linear mapping $\sigma$ from $\mathbf{Q}^{w}$ to $\mathbf{Q}^{w+1}$
$$ \sigma : x\mapsto x\ C.$$
The image space $\sigma(\mathbf{Q}^w)$ is a vector subspace of $\mathbf{Q}^{w+1}$. Let $\gamma_1, \gamma_2, \cdots, \gamma_w$ be the row vectors  of $C$. Thus the row space $R(C)$ is subspace of $\mathbf{Q}^{w+1}$ spanned  by $\gamma_1, \gamma_2, \cdots, \gamma_w$. So $\sigma(\mathbf{Q}^w)=R(C)$. By Remark 2.17, since  $C$ is of full row rank,  $dim_{\mathbf{Q}}(R(C))=w.$ So $\sigma$ is an one-one linear mapping  from $\mathbf{Q}^w$ to $R(C)$.

The equation (16) is an identity in $x_1, x_2, \cdots,x_{w} .$
Each of the f's is a rational combination of the x's, since $f=x\, C$.  By Remark 2.17, since  $C$ is of full row rank, each of the x's is a rational combination of the f's. Thus the equation (16) is an identity in
 the variables $f_1, f_2, \cdots,f_{w}, f_{w+1}$ for any  $f\in R(C)$.

 We express the integer $\alpha$ as the sum of two squares by the assumption, and bracket the terms $x_1^2+\cdots +x_{w}^2$  in twos.
 Each product of sums of two squares is itself a sum of two squares, i. e. Lemma 3.4, and  so (16) yields $$f_1^2+f_2^2+\cdots + f_{w}^2+f_{w+1}^2$$
\begin{equation}
 = y_1^2+y_2^2+\cdots+y_{w}^2+ \beta  z^2,
\end{equation}
where $z=x_1+x_2+\cdots+x_{w}$, and the $y's$ are related to the
$x's$ by an invertible linear transformation with rational
coefficients. Since the $y's$ are rational linear combinations of
the $x's$, it follows that the $y's $ (and $z$) are rational linear
combinations of the $f's$.  Thus the equation (17) is an identity in
 the variables $f_1, f_2, \cdots,f_{w}, f_{w+1}$ for any  $f\in R(C)$.

Suppose that  $ y_i=b_{i\, 1}f_1+\cdots +b_{i\, w}f_{w}+b_{i\, w+1}f_{w+1},\, 1 \leq i \leq w $.  We can
define $f_1$ as a rational linear combination of $f_2,\cdots ,
f_{w+1}$, in such a way that $y_1^2=f_1^2$: if $b_{1\, 1}\neq 1$ we
set $f_1=\frac{1}{1-b_{1\, 1}}(b_{1\, 2}f_2+\cdots+b_{1\,
w+1}f_{w+1})$, while if  $b_{1\, 1}= 1$ we set
$f_1=\frac{1}{-1-b_{1\, 1}}(b_{1\, 2}f_2+\cdots+b_{1\,
w+1}f_{w+1})$. Now we know that  $y_2$ is  a rational linear
combination of the $f's$, and, using the relevant expression for
$f_1$ found above, we can express $y_2$  as a rational linear
combination of $f_2,\cdots,f_{w+1}$. As before, we fix $f_2$ as a
rational combination of $f_3,\cdots , f_{w+1}$ in such a way that
$y_2^2=f_2^2$. Continuing thus, we eventually obtain $y_1,\cdots
,y_{w}$ and $f_1,\cdots ,f_{w}$ as rational multiples  of
$f_{w+1}$, satisfying $f_i^2=y_i^2\, (1 \leq i \leq w)$.

We reduce the equations step by step in this way until a truncated triangle of equations is obtained, say
 $$f_1=d_{1\, 2}f_2+\cdots+d_{1\,
w+1}f_{w+1},  $$
$$f_2=d_{2\, 3}f_3+\cdots+d_{2\,
w+1}f_{w+1},  $$
$$\cdots \cdots\cdots\cdots\cdots   $$
$$f_w=d_{w\, \, w+1}f_{w+1};  $$
$$f_i^2=y_i^2, \, (1 \leq i \leq w) ;$$
where $d_{i\, j} \in \mathbf{ Q}$.

For any  $x\in \mathbf{Q}^w$ and $x\neq 0$,  by Remark 2.17, since  $C$ is of full row rank, $f=x\, C$ implies $f\neq 0$. Let the last one $f_{w+1}\neq 0$ with suitable renumberings, if necessary.
Choose any non-zero rational value for  $f_{w+1} $.
All the
$y's $, the remaining $f's$, and $z$, are determined as above, and
substituting these values in (17) we obtain
\begin{equation}f_{w+1}^2= \beta z^2.
\end{equation}
Multiplying by a suitable constant  we have  that $\beta $ is a perfect squarer. So the theorem is proved.  \pfend

\begin{exa}
 There is a  symmetric $(45,12, 3)$ design (see\cite{l5}).
Let $A$ be its incidence matrix, which is  a 45 by 45 matrix.
Choose its  bordered  matrix $C$ is  a 46 by 47 matrix
as the following matrix.
 \[\begin{array}{ccccc} C=\left[\begin{array}{ccccc} A_{1\,\, 1} \ A_{1\,\, 2}
\\ A_{2\,\, 1}\ A_{2\,\, 2}
\end{array}\right ],
\end{array}\]
$$A_{1\,\, 1}=A .$$ $A_{1\,\, 2}$ is a 45 by 2 matrix  and
$$ A_{1\,\, 2}
=(1*1_{45}^t,  0*1_{45}^t). $$
 $A_{2\,\, 1}$ is
 a 1 by 45 matrix and
  \[\begin{array}{ccccc} A_{2\,\, 1}=\left[\begin{array}{ccccc}    \frac{2}{9}\cdot1_{45}
\end{array}\right ].
\end{array}\]
$A_{2\,\, 2}$ is a 1 by 2 matrix and
 \[ \begin{array}{ccccc} A_{2\,\, 2}=\left[\begin{array}{ccccccccccccccc}    \frac{4}{3} \ 3   \end{array}\right ] . \end{array}
\]
  It is easy to check that $C$ has  the property of row inner products, i.e.,

  (i)\,\, the inner product of any two distinct rows of $C$ is equal to $4$ ;

  (ii)\,\, and
  the inner product of any rows with themselves of $C$ is equal to $13$.\linebreak
  It follows that
 $$ C\, C^t=9 I_{46}+ 4 J_{46}$$
  and $C$ is exactly  the  bordered matrix of the symmetric $(45,12, 3)$ design.
 We have that 4  is  a  square   just as the assertion of the above theorem.
\end{exa}

\begin{lem}  (\textbf{Case 4 of Main Theorem 1} ) Let  $C$ be  a $w$ by $w+2$ nonsquare rational  matrix without any column of $k\cdot 1^t_{w+2}$, where $1^t_{w+2}$ is the $w+2$-dimensional all 1 column vector and $k$ is a rational number, $\alpha, \beta$ be positive  integers and  $\alpha=a^2+b^2$, where $a, b$ are  integers.
 Suppose the matrix $\alpha I_{w}+\beta J_{w}$ is the  positive definite matrix with plus 2 congruent factorization  property such that
\begin{equation} C\, C^t=\alpha I_{w}+\beta J_{w}.\end{equation}
If  $w\equiv 2\, (mod\, 4) $, then  $\beta $ is a sum of two squares.
\end{lem}
\pf
By the assumption   we have the identity
$$
 C \, \, C^t= \alpha I_{w}+  \beta  J_{w}$$ for the rational  matrix $C$.
 The idea of the proof is to interpret this as an identity in quadratic forms over the rational field.

  Suppose that   $w\equiv 2\, (mod\, 4)$.
  If $x$ is the row vector $(x_1,x_2,\cdots, x_{w})$, then the
  identity for $C\, C^t$ gives $$xC\, C^t x^t= \alpha(x_1^2+x_2^2+\cdots+x_{w}^2)+\beta(x_1+x_2+\cdots+x_{w})^2. $$
Putting $f=x\, C$, $f=(f_1,f_2,\cdots, f_{w}, f_{w+1},  f_{w+2})$, $z=x_1+x_2+\cdots+x_{w}$, we have $f\, f^t= x C\, C^tx^t$ and
\begin{equation} f_1^2+f_2^2+\cdots+f_{w}^2+f_{w+1}^2 +f_{w+2}^2
=
\alpha(x_1^2+x_2^2+\cdots+x_{w}^2)+\beta z^2;
\end{equation}
$$f_1=c_{1\, 1}x_1+c_{2\, 1}x_2\cdots+c_{w\, 1} x_{w},  $$
$$f_2=c_{1\, 2}x_1+c_{2\, 2}x_2\cdots+c_{w\, 2} x_{w},  $$
$$\cdots \cdots\cdots\cdots\cdots   $$
$$f_w=c_{1\, w}x_1+c_{2\, w}x_2\cdots+c_{w\, w} x_{w},  $$ $$f_{w+1}=c_{1\, {w+1}}x_1+c_{2\, {w+1}}x_2\cdots+c_{w\, {w+1}} x_{w},  $$
$$f_{w+2}=c_{1\, {w+2}}x_1+c_{2\, {w+2}}x_2\cdots+c_{w\, {w+2}} x_{w};  $$
$$z=x_1+x_2+\cdots+x_{w} .$$
Thus the cone (20) of variables $f_1, f_2, \cdots, f_{w}, f_{w+1},f_{w+2}, x_1,x_2, \cdots, x_{w}, z$ has some nontrivial rational points. We will get a nontrivial rational point for $f_{w+1}^2+f_{w+2}^2= \beta z^2$  such that $f_{w+2}^2 \neq 0$ by the Ryser-Chowla elimination procedure for the above homogeneous equations.

Now define a linear mapping $\sigma$ from $\mathbf{Q}^{w}$ to $\mathbf{Q}^{w+2}$
$$ \sigma : x\mapsto x\ C.$$
The image space $\sigma(\mathbf{Q}^w)$ is a vector subspace of $\mathbf{Q}^{w+2}$. Let $\gamma_1, \gamma_2, \cdots, \gamma_w$ be the row vectors  of $C$. Thus the row space $R(C)$ is subspace of $\mathbf{Q}^{w+2}$ spanned  by $\gamma_1, \gamma_2, \cdots, \gamma_w$. So $\sigma(\mathbf{Q}^w)=R(C)$. By Remark 2.17, since  $C$ is of full row rank,  $dim_{\mathbf{Q}}(R(C))=w.$ So $\sigma$ is an one-one linear mapping  from $\mathbf{Q}^w$ to $R(C)$.

The equation (20) is an identity in $x_1, x_2, \cdots,x_{w} .$
Each of the f's is a rational combination of the x's, since $f=x\, C$. By Remark 2.17, since  $C$ is of full row rank, each of the x's is a rational combination of the f's. Thus the equation (20) is an identity in
 the variables $f_1, f_2, \cdots,f_{w}, f_{w+1}, f_{w+2}$ for any  $f\in R(C)$.

 We express the integer $\alpha$ as the sum of two squares by the assumption, and bracket the terms $x_1^2+\cdots +x_{w}^2$  in twos.
 Each product of sums of two squares is itself a sum of two squares, and  so (20) yields $$f_1^2+f_2^2+\cdots + f_{w}^2+f_{w+1}^2+f_{w+2}^2$$
\begin{equation}
 = y_1^2+y_2^2+\cdots+y_{w}^2+ \beta  z^2,
\end{equation}
where $z=x_1+x_2+\cdots+x_{w}$, and the $y's$ are related to the
$x's$ by an invertible linear transformation with rational
coefficients. Since the $y's$ are rational linear combinations of
the $x's$, it follows that the $y's $ (and $z$) are rational linear
combinations of the $f's$.  Thus the equation (21) is an identity in
 the variables $f_1, f_2, \cdots,f_{w}, f_{w+1}, f_{w+2}$ for any  $f\in R(C)$.

Suppose that   $ y_i=b_{i\, 1}f_1+\cdots +b_{i\, w}f_{w}+b_{i\, w+1}f_{w+1} +b_{i\, w+2}f_{w+2}, \, 1 \leq i \leq w $.
We can
define $f_1$ as a rational linear combination of $f_2,\cdots ,
f_{w+1}, f_{w+2}$, in such a way that $y_1^2=f_1^2$: if $b_{1\, 1}\neq 1$ we
set $f_1=\frac{1}{1-b_{1\, 1}}(b_{1\, 2}f_2+\cdots+b_{1\,
w+1}f_{w+1 +b_{1\, w+2}f_{w+2}})$, while if  $b_{1\, 1}= 1$ we set
$f_1=\frac{1}{-1-b_{1\, 1}}(b_{1\, 2}f_2+\cdots+b_{1\,
w+1}f_{w+1} +b_{1\, w+2}f_{w+2})$. Now we know that  $y_2$ is  a rational linear
combination of the $f's$, and, using the relevant expression for
$f_1$ found above, we can express $y_2$  as a rational linear
combination of $f_2,\cdots,f_{w+1},f_{w+2}$. As before, we fix $f_2$ as a
rational combination of $f_3,\cdots , f_{w+1}, f_{w+2}$ in such a way that
$y_2^2=f_2^2$. Continuing thus, we eventually obtain $y_1,\cdots
,y_{w}$ and $f_1,\cdots ,f_{w}$ as rational combinations  of
$f_{w+1}, f_{w+2}$, satisfying $f_i^2=y_i^2\, (1 \leq i \leq w)$.

We reduce the equations step by step in this way until a truncated triangle of equations is obtained, say
 $$f_1=d_{1\, 2}f_2+\cdots+d_{1\,
w+1}f_{w+1}+d_{1\,
w+2}f_{w+2},  $$
$$f_2=d_{2\, 3}f_3+\cdots+d_{2\,
w+1}f_{w+1}+d_{2\,
w+2}f_{w+2},  $$
$$\cdots \cdots\cdots\cdots\cdots   $$
$$f_w=d_{w\, \, w+1}f_{w+1}+d_{w\,
w+2}f_{w+2}; $$
$$f_i^2=y_i^2, \, (1 \leq i \leq w); $$
where $d_{i\, j} \in \mathbf{ Q}$.

For any  $x\in \mathbf{Q}^w$ and $x\neq 0$,  by Remark 2.17, since  $C$ is of full row rank, $f=x\, C$ implies $f\neq 0$. Let the last one $f_{w+2}\neq 0$ with suitable renumberings, if necessary.
Choose any non-zero rational value for $f_{w+2} $.
All the
$y's $, the remaining $f's$, and $z$, are determined as above, and
substituting these values in (21) we obtain
\begin{equation}f_{w+1}^2+f_{w+2}^2= \beta z^2.
\end{equation}
Multiplying by a suitable constant  we have  that $\beta $ is a sum of two squares. So the theorem is proved.  \pfend
\begin{exa}
 There  is a symmetric $(45,12, 3)$ design (see\cite{l5}).
Let $A$ be its incidence matrix, which is  a 45 by 45 matrix.
Choose its  bordered matrix $C$ is  a 46 by 48 matrix
as the following matrix.
 \[\begin{array}{ccccc} C=\left[\begin{array}{ccccc} A_{1\,\, 1} \ A_{1\,\, 2}
\\ A_{2\,\, 1}\ A_{2\,\, 2}
\end{array}\right ],
\end{array}\]
$$A_{1\,\, 1}=A .$$ $A_{1\,\, 2}$ is a 45 by 3 matrix  and
$$ A_{1\,\, 2}
=(1*1_{45}^t,  1*1_{45}^t,0*1_{45}^t). $$
 $A_{2\,\, 1}$ is
 a 1 by 45 matrix and
  \[\begin{array}{ccccc} A_{2\,\, 1}=\left[\begin{array}{ccccc}    0 \cdot1_{45}
\end{array}\right ].
\end{array}\]
$A_{2\,\, 2}$ is a 1 by 3 matrix and
 \[ \begin{array}{ccccc} A_{2\,\, 2}=\left[\begin{array}{ccccccccccccccc}  3 \  2 \ 1   \end{array}\right ] . \end{array}
\]
  It is easy to check that $C$ has  the property of row inner products, i.e.,

  (i)\,\, the inner product of any two distinct rows of $C$ is equal to $5$ ;

  (ii)\,\, and
  the inner product of any rows with themselves of $C$ is equal to $14$.\linebreak
  It follows that
 $$ C\, C^t=9 I_{46}+ 5 J_{46}$$
  and $C$ is exactly  the  bordered matrix of the symmetric $(45,12, 3)$ design.
 We have that 5  is  is a sum of two squares  just as the assertion of the above theorem.
\end{exa}

\begin{lem}  (\textbf{Case 5 of Main Theorem 1} ) Let  $C$ be  a $w$ by $w+1$ nonsquare rational  matrix without any column of $k\cdot 1^t_{w+1}$, where $1^t_{w+1}$ is the $w+1$-dimensional all 1 column vector and $k$ is a rational number, $\alpha, \beta$ be positive  integers and
 $\alpha=a^2+b^2$, where $a, b$ are   integers.
 Suppose the matrix $\alpha I_{w}+\beta J_{w}$ is the  positive definite matrix with plus 1 congruent factorization  property such that
 \begin{equation} C\, C^t=\alpha I_{w}+\beta J_{w}.\end{equation}
If  $w\equiv 1\, (mod\, 4) $, then  $\alpha^*=\beta^* $.
\end{lem}
\pf By the assumption   we have the identity
$$
 C \, \, C^t= \alpha I_{w}+  \beta  J_{w}$$ for the rational  matrix $C$.
 The idea of the proof is to interpret this as an identity in quadratic forms over the rational field.

  Suppose that   $w\equiv 1\, (mod\, 4)$.
  If $x$ is the row vector $(x_1,x_2,\cdots, x_{w})$, then the
  identity for $C\, C^t$ gives $$xC\, C^t x^t= \alpha(x_1^2+x_2^2+\cdots+x_{w}^2)+\beta(x_1+x_2+\cdots+x_{w})^2. $$
Putting $f=x\, C$, $f=(f_1,f_2,\cdots, f_{w}, f_{w+1})$, $z=x_1+x_2+\cdots+x_{w}$, we have $f\, f^t= x C\, C^tx^t$ and
\begin{equation}  f_1^2+f_2^2+\cdots+f_{w}^2+f_{w+1}^2
=
\alpha(x_1^2+x_2^2+\cdots+x_{w}^2)+\beta z^2;
\end{equation}
$$f_1=c_{1\, 1}x_1+c_{2\, 1}x_2\cdots+c_{w\, 1} x_{w},  $$
$$f_2=c_{1\, 2}x_1+c_{2\, 2}x_2\cdots+c_{w\, 2} x_{w},  $$
$$\cdots \cdots\cdots\cdots\cdots   $$
$$f_w=c_{1\, w}x_1+c_{2\, w}x_2\cdots+c_{w\, w} x_{w},  $$
$$f_{w+1}=c_{1\, {w+1}}x_1+c_{2\, {w+1}}x_2\cdots+c_{w\, {w+1}} x_{w};  $$
$$z=x_1+x_2+\cdots+x_{w} .$$
Thus the cone (24) of variables $f_1, f_2, \cdots, f_{w}, f_{w+1}, x_1,x_2, \cdots, x_{w}, z$ has some nontrivial rational points. We will get a nontrivial rational point for $\alpha y^2=\beta z^2$ such that $y\neq 0$ by the Ryser-Chowla elimination procedure for the above homogeneous equations.

Now define a linear mapping $\sigma$ from $\mathbf{Q}^{w}$ to $\mathbf{Q}^{w+1}$
$$ \sigma : x\mapsto x\ C.$$
The image space $\sigma(\mathbf{Q}^w)$ is a vector subspace of $\mathbf{Q}^{w+1}$. Let $\gamma_1, \gamma_2, \cdots, \gamma_w$ be the row vectors  of $C$. Thus the row space $R(C)$ is subspace of $\mathbf{Q}^{w+1}$ spanned  by $\gamma_1, \gamma_2, \cdots, \gamma_w$. So $\sigma(\mathbf{Q}^w)=R(C)$. By Remark 2.17, since  $C$ is of full row rank,  $dim_{\mathbf{Q}}(R(C))=w.$ So $\sigma$ is an one-one linear mapping  from $\mathbf{Q}^w$ to $R(C)$.

The equation (24) is an identity in $x_1, x_2, \cdots,x_{w} .$
Each of the f's is a rational combination of the x's, since $f=x\, C$. By Remark 2.17, since  $C$ is of full row rank, each of the x's is a rational combination of the f's. Thus the equation (24) is an identity in
 the variables $f_1, f_2, \cdots,f_{w}, f_{w+1}$ for any  $f\in R(C)$.

 We express the integer $\alpha$ as the sum of two squares by the assumption, and bracket the terms $f_1^2+\cdots +f_{w}^2+f_{w+1}^2$  in twos.
 Each product of sums of two squares is itself a sum of two squares, and  so (24) yields
  $$\alpha(y_1^2+y_2^2+\cdots+y_{w}^2+y_{w+1}^2)  $$
\begin{equation}
=
\alpha(x_1^2+x_2^2+\cdots+x_{w}^2)+\beta z^2,
\end{equation}
 where $z=x_1+x_2+\cdots+x_{w}$, and the $y's$ are related to the
$f's$ by an invertible linear transformation with rational
coefficients. Thus $det(P)\neq 0$, $y=f \, P$.

Now define a linear mapping $\tau$ from $\mathbf{Q}^{w+1}$ to $\mathbf{Q}^{w+1}$
$$ \tau  : f\mapsto f\ P.$$
The image space $\tau \sigma(\mathbf{Q}^w)=V$ is a vector subspace of $\mathbf{Q}^{w+1}$ and $dim_{\mathbf{Q}}V=w$.

Since the $x's$ are rational linear combinations of
the $f's$, it follows that the $x's $ (and $z$) are rational linear
combinations of the $y's$. 
Thus the equation (25) is an identity in
 the variables $y_1, y_2, \cdots,y_{w}, y_{w+1}$ for any  $y\in V $.

Suppose that $ x_i=b_{i\, 1}y_1+\cdots +b_{i\, w}y_{w}+b_{i\, w+1}y_{w+1}$, $1\leq i \leq w$. We can
define $y_1$ as a rational linear combination of $y_2,\cdots ,
y_{w+1}$, in such a way that $x_1^2=y_1^2$: if $b_{1\, 1}\neq 1$ we
set $y_1=\frac{1}{1-b_{1\, 1}}(b_{1\, 2}y_2+\cdots+b_{1\,
w+1}y_{w+1})$, while if  $b_{1\, 1}= 1$ we set
$y_1=\frac{1}{-1-b_{1\, 1}}(b_{1\, 2}y_2+\cdots+b_{1\,
w+1}y_{w+1})$. Now we know that  $x_2$ is  a rational linear
combination of the $y's$, and, using the relevant expression for
$y_1$ found above, we can express $x_2$  as a rational linear
combination of $y_2,\cdots,y_{w+1}$. As before, we fix $x_2$ as a
rational combination of $y_3,\cdots , y_{w+1}$ in such a way that
$x_2^2=y_2^2$. Continuing thus, we eventually obtain $x_1,\cdots
,x_{w}$ and $y_1,\cdots ,y_{w}$ as rational multiples  of
$y_{w+1}$, satisfying $x_i^2=y_i^2\, (1 \leq i \leq w)$.

We reduce the equations step by step in this way until a truncated triangle of equations is obtained, say
 $$y_1=d_{1\, 2}y_2+\cdots+d_{1\,
w+1}y_{w+1},  $$
$$y_2=d_{2\, 3}y_3+\cdots+d_{2\,
w+1}y_{w+1},  $$
$$\cdots \cdots\cdots\cdots\cdots   $$
$$y_w=d_{w\, \, w+1}y_{w+1};  $$
$$x_i^2=y_i^2, \, (1 \leq i \leq w); $$
where $d_{i\, j} \in \mathbf{ Q}$.

For any  $x\in \mathbf{Q}^w$ and $x\neq 0$,  by Remark 2.17, since  $C$ is of full row rank, $f=x\, C$, and $det(P)\neq 0$, $y=f \, P$, it implies $y\neq 0$. Let the last one $y_{w+1}\neq 0$ with suitable renumberings, if necessary.
Choose any non-zero rational value for  $y_{w+1} $.
All the
$x's $, the remaining $y's$, and $z$, are determined as above, and
substituting these values in (25) we obtain
\begin{equation}\alpha y_{w+1}^2= \beta z^2.
\end{equation}
Multiplying by a suitable constant  we have  that  $\alpha^*=\beta^* $. So the theorem is proved.  \pfend
\begin{exa}
 There  is a symmetric $(36,15, 6)$ design (see\cite{l5}).
Let $A$ be its incidence matrix, which is  a 36 by 36 matrix.
Choose its  bordered matrix $C$ is  a 37 by 38 matrix
as the following matrix.
 \[\begin{array}{ccccc} C=\left[\begin{array}{ccccc} A_{1\,\, 1} \ A_{1\,\, 2}
\\ A_{2\,\, 1}\ A_{2\,\, 2}
\end{array}\right ],
\end{array}\]
$$A_{1\,\, 1}=A .$$ $A_{1\,\, 2}$ is a 36 by 2 matrix  and
$$ A_{1\,\, 2}
=(3*1_{36}^t,  1*1_{36}^t). $$
 $A_{2\,\, 1}$ is
 a 1 by 36 matrix and
  \[\begin{array}{ccccc} A_{2\,\, 1}=\left[\begin{array}{ccccc}    \frac{7}{9} \cdot1_{36}
\end{array}\right ].
\end{array}\]
$A_{2\,\, 2}$ is a 1 by 2 matrix and
 \[ \begin{array}{ccccc} A_{2\,\, 2}=\left[\begin{array}{ccccccccccccccc}  \frac{14}{15} \  \frac{23}{15}   \end{array}\right ] . \end{array}
\]
  It is easy to check that $C$ has  the property of row inner products, i.e.,

  (i)\,\, the inner product of any two distinct rows of $C$ is equal to $16$ ;

  (ii)\,\, and
  the inner product of any rows with themselves of $C$ is equal to $25$.\linebreak
  It follows that
 $$ C\, C^t=9 I_{37}+ 16 J_{37}$$
  and $C$ is exactly  the  bordered matrix of the symmetric $(36,15, 6)$ design.
 We have that $9^*= 16^*$  just as the assertion of the above theorem.
\end{exa}
\begin{lem}  (\textbf{Case 6 of Main Theorem 1} ) Let  $C$ be  a $w$  by $w+2$ nonsquare rational  matrix without any column of $k\cdot 1^t_{w+2}$, where $1^t_{w+2}$ is the $w+2$-dimensional all 1 column vector and $k$ is a rational number,   $\alpha, \beta$ be positive  integers and
  $\alpha=a^2+b^2$, where $a, b$ are   integers.
 Suppose  matrix $\alpha I_{w}+\beta J_{w}$ is the  positive definite matrix with plus $2$  congruent factorization property such that
 \begin{equation} C\, C^t=\alpha I_{w}+\beta J_{w}.\end{equation}
If  $w\equiv 1\, (mod\, 4) $, then    the equation
  $$ \alpha z^2= -x^2+ \beta y^2$$
must  have a solution in integers,
$x, y, z$,  not all zero.
\end{lem}
\pf By the assumption   we have the identity
$$
 C \, \, C^t= \alpha I_{w}+  \beta  J_{w}$$ for the rational  matrix $C$.
 The idea of the proof is to interpret this as an identity in quadratic forms over the rational field.

  Suppose that   $w\equiv 1\, (mod\, 4)$.
  If $x$ is the row vector $(x_1,x_2,\cdots, x_{w})$, then the
  identity for $C\, C^t$ gives $$xC\, C^t x^t= \alpha(x_1^2+x_2^2+\cdots+x_{w}^2)+\beta(x_1+x_2+\cdots+x_{w})^2. $$
Putting $f=x\, C$,  $f=(f_1,f_2,\cdots, f_{w}, f_{w+1}, f_{w+2})$, $z=x_1+x_2+\cdots+x_{w}$, we have $f\, f^t= x C\, C^tx^t$ and
\begin{equation}f_1^2+f_2^2+\cdots+f_{w}^2+f_{w+1}^2 +f_{w+2}^2
=
\alpha(x_1^2+x_2^2+\cdots+x_{w}^2)+\beta z^2;
\end{equation}
$$f_1=c_{1\, 1}x_1+c_{2\, 1}x_2\cdots+c_{w\, 1} x_{w},  $$
$$f_2=c_{1\, 2}x_1+c_{2\, 2}x_2\cdots+c_{w\, 2} x_{w},  $$
$$\cdots \cdots\cdots\cdots\cdots   $$
$$f_w=c_{1\, w}x_1+c_{2\, w}x_2\cdots+c_{w\, w} x_{w},  $$
$$f_{w+1}=c_{1\, {w+1}}x_1+c_{2\, {w+1}}x_2\cdots+c_{w\, {w+1}} x_{w}, $$
$$f_{w+2}=c_{1\, {w+2}}x_1+c_{2\, {w+2}}x_2\cdots+c_{w\, {w+2}} x_{w};  $$
$$z=x_1+x_2+\cdots+x_{w} .$$
Thus the cone (28) of variables $f_1, f_2, \cdots, f_{w}, f_{w+1}, f_{w+2}, x_1,x_2, \cdots, x_{w}, z$ has some nontrivial rational points. We will get a nontrivial rational point for $\alpha y^2+x^2=\beta z^2$ such that $x\neq 0$ by the Ryser-Chowla elimination procedure for the above homogeneous equations.

Now define a linear mapping $\sigma$ from $\mathbf{Q}^{w}$ to $\mathbf{Q}^{w+2}$
$$ \sigma : x\mapsto x\ C.$$
The image space $\sigma(\mathbf{Q}^w)$ is a vector subspace of $\mathbf{Q}^{w+2}$. Let $\gamma_1, \gamma_2, \cdots, \gamma_w$ be the row vectors  of $C$. Thus the row space $R(C)$ is subspace of $\mathbf{Q}^{w+2}$ spanned  by $\gamma_1, \gamma_2, \cdots, \gamma_w$. So $\sigma(\mathbf{Q}^w)=R(C)$. By Remark 2.17, since  $C$ is of full row rank,  $dim_{\mathbf{Q}}(R(C))=w.$ So $\sigma$ is an one-one linear mapping  from $\mathbf{Q}^w$ to $R(C)$.

The equation (28) is an identity in $x_1, x_2, \cdots,x_{w} .$
Each of the f's is a rational combination of the x's, since $f=x\, C$. By Remark 2.17, since  $C$ is of full row rank, each of the x's is a rational combination of the f's. Thus the equation (28) is an identity in
 the variables $f_1, f_2, \cdots,f_{w}, f_{w+1}, f_{w+2}$ for any  $f\in R(C)$.

 We express the integer $\alpha$ as the sum of two squares by the assumption, and bracket the terms $f_1^2+\cdots +f_{w}^2+f_{w+1}^2$  in twos.
 Each product of sums of two squares is itself a sum of two squares, and  so (28) yields
  $$\alpha(y_1^2+y_2^2+\cdots+y_{w}^2+y_{w+1}^2)+y_{w+2}^2  $$
\begin{equation}
=
\alpha(x_1^2+x_2^2+\cdots+x_{w}^2)+\beta z^2,
\end{equation}
 where $z=x_1+x_2+\cdots+x_{w}$, $f_{w+2}=y_{w+2}$, and the $y's$ are related to the
$f's$ by an invertible linear transformation with rational
coefficients.  Thus $det(P)\neq 0$, $y=f \, P$.

Now define a linear mapping $\tau$ from $\mathbf{Q}^{w+2}$ to $\mathbf{Q}^{w+2}$
$$ \tau  : f\mapsto f\ P.$$
The image space $\tau \sigma(\mathbf{Q}^w)=V$ is a vector subspace of $\mathbf{Q}^{w+2}$ and $dim_{\mathbf{Q}}V=w$.

Since the $x's$ are rational linear combinations of
the $f's$, it follows that the $x's $ (and $z$) are rational linear
combinations of the $y's$.  Thus the equation (29) is an identity in
 the variables $y_1, y_2, \cdots,y_{w}, y_{w+1},y_{w+2}$  for any  $y\in V $.

Suppose that $ x_i=b_{i\, 1}y_1+\cdots +b_{i\, w}y_{w}+b_{i\, w+1}y_{w+1}+b_{i\, w+2}y_{w+2}$, $1\leq i\leq w$. We can
define $y_1$ as a rational linear combination of $y_2,\cdots ,
y_{w+1},y_{w+2}$, in such a way that $x_1^2=y_1^2$: if $b_{1\, 1}\neq 1$ we
set $y_1=\frac{1}{1-b_{1\, 1}}(b_{1\, 2}y_2+\cdots+b_{1\,
w+1}y_{w+1}+b_{1\, w+2}y_{w+2})$, while if  $b_{1\, 1}= 1$ we set
$y_1=\frac{1}{-1-b_{1\, 1}}(b_{1\, 2}y_2+\cdots+b_{1\,
w+1}y_{w+1}+b_{1\, w+2}y_{w+2})$. Now we know that  $x_2$ is  a rational linear
combination of the $y's$, and, using the relevant expression for
$y_1$ found above, we can express $x_2$  as a rational linear
combination of $y_2,\cdots,y_{w+1},y_{w+2}$. As before, we fix $x_2$ as a
rational combination of $y_3,\cdots , y_{w+1},y_{w+2}$ in such a way that
$x_2^2=y_2^2$. Continuing thus, we eventually obtain $x_1,\cdots
,x_{w}$ and $y_1,\cdots ,y_{w}$ as rational linear combinations   of
$y_{w+1},y_{w+2}$, satisfying $x_i^2=y_i^2\, (1 \leq i \leq w)$.

We reduce the equations step by step in this way until a truncated triangle of equations is obtained, say
 $$y_1=d_{1\, 2}y_2+\cdots+d_{1\,
w+1}y_{w+1}+d_{1\,
w+2}y_{w+2},  $$
$$y_2=d_{2\, 3}y_3+\cdots+d_{2\,
w+1}y_{w+1}+d_{2\,
w+2}y_{w+2},  $$
$$\cdots \cdots\cdots\cdots\cdots   $$
$$y_w=d_{w\, \, w+1}y_{w+1}+d_{w\, \, w+2}y_{w+2};  $$
$$x_i^2=y_i^2, \, (1 \leq i \leq w) ;$$
where $d_{i\, j} \in \mathbf{ Q}$.

For any  $x\in \mathbf{Q}^w$ and $x\neq 0$,  by Remark 2.17, since  $C$ is of full row rank, $f=x\, C$, and $det(P)\neq 0$, $y=f \, P$, it implies $y\neq 0$. Let the last one $y_{w+2}\neq 0$ with suitable renumberings, if necessary.
Choose any non-zero rational value for  $y_{w+2} $.
All the
$x's $, the remaining $y's$, and $z$, are determined as above, and
substituting these values in (29) we obtain
\begin{equation}\alpha y_{w+1}^2+y_{w+2}^2= \beta z^2.
\end{equation}
Multiplying by a suitable constant  we have  that   $$ \alpha z^2= -x^2+ \beta y^2.$$  So the theorem is proved.  \pfend
\begin{exa}
 There  is a symmetric $(36,15, 6)$ design (see\cite{l5}).
Let $A$ be its incidence matrix, which is  a 36 by 36 matrix.
Choose its  bordered matrix $C$ is  a 37 by 39 matrix
as the following matrix.
 \[\begin{array}{ccccc} C=\left[\begin{array}{ccccc} A_{1\,\, 1} \ A_{1\,\, 2}
\\ A_{2\,\, 1}\ A_{2\,\, 2}
\end{array}\right ],
\end{array}\]
$$A_{1\,\, 1}=A .$$ $A_{1\,\, 2}$ is a 36 by 3 matrix  and
$$ A_{1\,\, 2}
=(2*1_{36}^t,  0*1_{36}^t,   0*1_{36}^t). $$
 $A_{2\,\, 1}$ is
 a 1 by 36 matrix and
  \[\begin{array}{ccccc} A_{2\,\, 1}=\left[\begin{array}{ccccc}    \frac{8}{15} \cdot1_{36}
\end{array}\right ].
\end{array}\]
$A_{2\,\, 2}$ is a 1 by 3 matrix and
 \[ \begin{array}{ccccc} A_{2\,\, 2}=\left[\begin{array}{ccccccccccccccc} 1\ 1  \  \frac{13}{5}   \end{array}\right ] . \end{array}
\]
  It is easy to check that $C$ has  the property of row inner products, i.e.,

  (i)\,\, the inner product of any two distinct rows of $C$ is equal to $10$ ;

  (ii)\,\, and
  the inner product of any rows with themselves of $C$ is equal to $19$.\linebreak
  It follows that
 $$ C\, C^t=9 I_{37}+10 J_{37}$$
  and $C$ is exactly  the  bordered matrix of the symmetric $(36,15, 6)$ design.
 We have that  the equation
  $$ 9  z^2= -x^2+ 10 y^2$$
must  have a solution in integers,
$x, y, z$,  not all zero  just as the assertion of the above theorem.
\end{exa}
\begin{lem} (\textbf{Case 7 of  Main Theorem 1} )  Let  $C$ be  a $w$ by $w+1$ nonsquare rational  matrix without any column of $k\cdot 1^t_{w+1}$, where $1^t_{w+1}$ is the $w+1$-dimensional all 1 column vector and $k$ is a rational number, $\alpha, \beta$ be positive  integers.  Suppose
  the matrix $\alpha I_{w}+\beta J_{w}$ is the  positive definite matrix with plus $1$  congruent factorization property such that
 \begin{equation} C\, C^t=\alpha I_{w}+\beta J_{w}.\end{equation}
If  $w\equiv 3\, (mod\, 4) $,  then $\alpha^*=\beta^*$.
\end{lem}
\pf By the assumption   we have the identity
$$
 C \, \, C^t= \alpha I_{w}+  \beta  J_{w}$$ for the rational  matrix $C$.
 The idea of the proof is to interpret this as an identity in quadratic forms over the rational field.

  Suppose that   $w\equiv 3\, (mod\, 4)$.
  If $x$ is the row vector $(x_1,x_2,\cdots, x_{w})$, then the
  identity for $C\, C^t$ gives $$xC\, C^t x^t= \alpha(x_1^2+x_2^2+\cdots+x_{w}^2)+\beta(x_1+x_2+\cdots+x_{w})^2. $$
Putting $f=x\, C$, $f=(f_1,f_2,\cdots, f_{w}, f_{w+1})$, $z=x_1+x_2+\cdots+x_{w}$, we have $f\, f^t= x C\, C^tx^t$ and
\begin{equation}f_1^2+f_2^2+\cdots+f_{w}^2+f_{w+1}^2
=
\alpha(x_1^2+x_2^2+\cdots+x_{w}^2)+\beta z^2;
\end{equation}
$$f_1=c_{1\, 1}x_1+c_{2\, 1}x_2\cdots+c_{w\, 1} x_{w},  $$
$$f_2=c_{1\, 2}x_1+c_{2\, 2}x_2\cdots+c_{w\, 2} x_{w},  $$
$$\cdots \cdots\cdots\cdots\cdots   $$
$$f_w=c_{1\, w}x_1+c_{2\, w}x_2\cdots+c_{w\, w} x_{w},  $$
$$f_{w+1}=c_{1\, {w+1}}x_1+c_{2\, {w+1}}x_2\cdots+c_{w\, {w+1}} x_{w};  $$
$$z=x_1+x_2+\cdots+x_{w} .$$
Thus the cone (32) of variables $f_1, f_2, \cdots, f_{w}, f_{w+1}, x_1,x_2, \cdots, x_{w}, z$ has some nontrivial rational points. We will get a nontrivial rational point for $\alpha y^2=\beta z^2$ such that $y\neq 0$ by the Ryser-Chowla elimination procedure for the above homogeneous equations.

Now define a linear mapping $\sigma$ from $\mathbf{Q}^{w}$ to $\mathbf{Q}^{w+1}$
$$ \sigma : x\mapsto x\ C.$$
The image space $\sigma(\mathbf{Q}^w)$ is a vector subspace of $\mathbf{Q}^{w+1}$. Let $\gamma_1, \gamma_2, \cdots, \gamma_w$ be the row vectors  of $C$. Thus the row space $R(C)$ is subspace of $\mathbf{Q}^{w+1}$ spanned  by $\gamma_1, \gamma_2, \cdots, \gamma_w$. So $\sigma(\mathbf{Q}^w)=R(C)$. By Remark 2.17, since  $C$ is of full row rank,  $dim_{\mathbf{Q}}(R(C))=w.$ So $\sigma$ is an one-one linear mapping  from $\mathbf{Q}^w$ to $R(C)$.

The equation (32) is an identity in $x_1, x_2, \cdots,x_{w} .$
Each of the f's is a rational combination of the x's, since $f=x\, C$. By Remark 2.17, since  $C$ is of full row rank, each of the x's is a rational combination of the f's. Thus the equation (32) is an identity in
 the variables $f_1, f_2, \cdots,f_{w}, f_{w+1}$ for any  $f\in R(C)$.

We express the integer $\alpha$ as the sum of four squares by Lemma 3.1, and bracket the terms $f_1^2+\cdots +f_{w}^2+f_{w+1}^2$  in fours.
 Each product of sums of four squares is itself a sum of four squares, and  so (32) yields
  $$\alpha(y_1^2+y_2^2+\cdots+y_{w}^2+y_{w+1}^2)  $$
\begin{equation}
=
\alpha(x_1^2+x_2^2+\cdots+x_{w}^2)+\beta z^2,
\end{equation}
 where $z=x_1+x_2+\cdots+x_{w}$, and the $y's$ are related to the
$f's$ by an invertible linear transformation with rational
coefficients. Thus $det(P)\neq 0$, $y=f \, P$.

Now define a linear mapping $\tau$ from $\mathbf{Q}^{w+1}$ to $\mathbf{Q}^{w+1}$
$$ \tau  : f\mapsto f\ P.$$
The image space $\tau \sigma(\mathbf{Q}^w)=V$ is a vector subspace of $\mathbf{Q}^{w+1}$ and $dim_{\mathbf{Q}}V=w$.

Since the $x's$ are rational linear combinations of
the $f's$, it follows that the $x's $ (and $z$) are rational linear
combinations of the $y's$. 
Thus the equation (33) is an identity in
 the variables $y_1, y_2, \cdots,y_{w}, y_{w+1}$ for any  $y\in V $

Suppose that $ x_i=b_{i\, 1}y_1+\cdots +b_{i\, w}y_{w}+b_{i\, w+1}y_{w+1}$, $1\leq i\leq w$. We can
define $y_1$ as a rational linear combination of $y_2,\cdots ,
y_{w+1}$, in such a way that $x_1^2=y_1^2$: if $b_{1\, 1}\neq 1$ we
set $y_1=\frac{1}{1-b_{1\, 1}}(b_{1\, 2}y_2+\cdots+b_{1\,
w+1}y_{w+1})$, while if  $b_{1\, 1}= 1$ we set
$y_1=\frac{1}{-1-b_{1\, 1}}(b_{1\, 2}y_2+\cdots+b_{1\,
w+1}y_{w+1})$. Now we know that  $x_2$ is  a rational linear
combination of the $y's$, and, using the relevant expression for
$y_1$ found above, we can express $x_2$  as a rational linear
combination of $y_2,\cdots,y_{w+1}$. As before, we fix $x_2$ as a
rational combination of $y_3,\cdots , y_{w+1}$ in such a way that
$x_2^2=y_2^2$. Continuing thus, we eventually obtain $x_1,\cdots
,x_{w}$ and $y_1,\cdots ,y_{w}$ as rational multiples  of
$y_{w+1}$, satisfying $x_i^2=y_i^2\, (1 \leq i \leq w)$.

We reduce the equations step by step in this way until a truncated triangle of equations is obtained, say
 $$y_1=d_{1\, 2}y_2+\cdots+d_{1\,
w+1}y_{w+1},  $$
$$y_2=d_{2\, 3}y_3+\cdots+d_{2\,
w+1}y_{w+1},  $$
$$\cdots \cdots\cdots\cdots\cdots   $$
$$y_w=d_{w\, \, w+1}y_{w+1};  $$
$$x_i^2=y_i^2, \, (1 \leq i \leq w); $$
where $d_{i\, j} \in \mathbf{ Q}$.

For any  $x\in \mathbf{Q}^w$ and $x\neq 0$,  by Remark 2.17, since  $C$ is of full row rank, $f=x\, C$, and $det(P)\neq 0$, $y=f \, P$, it implies $y\neq 0$. Let the last one $y_{w+1}\neq 0$ with suitable renumberings, if necessary.
Choose any non-zero rational value for  $y_{w+1} $.
All the
$x's $, the remaining $y's$, and $z$, are determined as above, and
substituting these values in (33) we obtain
\begin{equation}\alpha y_{w+1}^2= \beta z^2.
\end{equation}
Multiplying by a suitable constant  we have  that  $\alpha^*=\beta^* $. So the theorem is proved.  \pfend
\begin{exa}  The projective plane of order 7 is the symmetric $(57,8, 1)$ design.
Let $A$ be its incidence matrix, which is  a 57 by 57 matrix.
Choose its  bordered matrix $C$ is  a 59 by 60 matrix
as the following matrix.
 \[\begin{array}{ccccc} C=\left[\begin{array}{ccccc} A_{1\,\, 1} \ A_{1\,\, 2}
\\ A_{2\,\, 1}\ A_{2\,\, 2}
\end{array}\right ],
\end{array}\]
$$A_{1\,\, 1}=A .$$ $A_{1\,\, 2}$ is a 57 by 3 matrix  and
$$ A_{1\,\, 2}
=(2*1_{57}^t,  1*1_{57}^t,   1*1_{57}^t). $$
 $A_{2\,\, 1}$ is
 a 2 by 57 matrix and
  \[\begin{array}{ccccc} A_{2\,\, 1}=\left[\begin{array}{ccccc}    0 \cdot1_{57} \\  0 \cdot1_{57}
\end{array}\right ].
\end{array}\]
$A_{2\,\, 2}$ is a 2 by 3 matrix and
 \[ \begin{array}{ccccc} A_{2\,\, 2}=\left[\begin{array}{ccccccccccccccc} 1\ 3  \  2 \\ \frac{11}{5}\ -\frac{2}{5}  \  3  \end{array}\right ] . \end{array}
\]
  It is easy to check that $C$ has  the property of row inner products, i.e.,

  (i)\,\, the inner product of any two distinct rows of $C$ is equal to $7$ ;

  (ii)\,\, and
  the inner product of any rows with themselves of $C$ is equal to $14$.\linebreak
  It follows that
 $$ C\, C^t=7 I_{59}+7 J_{59}$$
  and $C$ is exactly  the  bordered matrix of the symmetric $(57,8, 1)$ design.
 We have that $\alpha=7, \beta=7$ and $\alpha^*=\beta^*$  just as the assertion of the above theorem.
 \end{exa}
\begin{lem} (\textbf{Case 8 of Main Theorem 1} )  Let  $C$ be  a $w$ by $w+2$ nonsquare rational  matrix without any column of $k\cdot 1^t_{w+2}$, where $1^t_{w+2}$ is the $w+2$-dimensional all 1 column vector and $k$ is a rational number, $\alpha, \beta$ be two positive  integers.
 Suppose the  matrix $\alpha I_{w}+\beta J_{w}$ is the  positive definite matrix with plus $2$  congruent factorization property such that
 \begin{equation} C\, C^t=\alpha I_{w}+\beta J_{w}.\end{equation}
If  $w\equiv 3\, (mod\, 4) $, then   the equation
  $$ \alpha z^2=- x^2+ \beta y^2$$
must  have a solution in integers,
$x, y, z$,  not all zero.
\end{lem}

\pf By the assumption   we have the identity
$$
 C \, \, C^t= \alpha I_{w}+  \beta  J_{w}$$ for the rational  matrix $C$.
 The idea of the proof is to interpret this as an identity in quadratic forms over the rational field.

  Suppose that   $w\equiv 3\, (mod\, 4)$.
  If $x$ is the row vector $(x_1,x_2,\cdots, x_{w})$, then the
  identity for $C\, C^t$ gives $$xC\, C^t x^t= \alpha(x_1^2+x_2^2+\cdots+x_{w}^2)+\beta(x_1+x_2+\cdots+x_{w})^2. $$
Putting $f=x\, C$, $f=(f_1,f_2,\cdots, f_{w}, f_{w+1}, f_{w+2})$, $z=x_1+x_2+\cdots+x_{w}$, we have $f\, f^t= x C\, C^tx^t$ and
\begin{equation} f_1^2+f_2^2+\cdots+f_{w}^2+f_{w+1}^2 +f_{w+2}^2
=
\alpha(x_1^2+x_2^2+\cdots+x_{w}^2)+\beta z^2;
\end{equation}
$$f_1=c_{1\, 1}x_1+c_{2\, 1}x_2\cdots+c_{w\, 1} x_{w},  $$
$$f_2=c_{1\, 2}x_1+c_{2\, 2}x_2\cdots+c_{w\, 2} x_{w},  $$
$$\cdots \cdots\cdots\cdots\cdots   $$
$$f_w=c_{1\, w}x_1+c_{2\, w}x_2\cdots+c_{w\, w} x_{w},  $$
$$f_{w+1}=c_{1\, {w+1}}x_1+c_{2\, {w+1}}x_2\cdots+c_{w\, {w+1}} x_{w}, $$
$$f_{w+2}=c_{1\, {w+2}}x_1+c_{2\, {w+2}}x_2\cdots+c_{w\, {w+2}} x_{w};  $$
$$z=x_1+x_2+\cdots+x_{w} .$$
Thus the cone (36) of variables $f_1, f_2, \cdots, f_{w}, f_{w+1}, f_{w+2}, x_1,x_2, \cdots, x_{w}, z$ has some nontrivial rational points. We will get a nontrivial rational point for $\alpha y^2+x^2=\beta z^2$ such that $x\neq 0$ by the Ryser-Chowla elimination procedure for the above homogeneous equations.

Now define a linear mapping $\sigma$ from $\mathbf{Q}^{w}$ to $\mathbf{Q}^{w+2}$
$$ \sigma : x\mapsto x\ C.$$
The image space $\sigma(\mathbf{Q}^w)$ is a vector subspace of $\mathbf{Q}^{w+2}$. Let $\gamma_1, \gamma_2, \cdots, \gamma_w$ be the row vectors  of $C$. Thus the row space $R(C)$ is subspace of $\mathbf{Q}^{w+2}$ spanned  by $\gamma_1, \gamma_2, \cdots, \gamma_w$. So $\sigma(\mathbf{Q}^w)=R(C)$. By Remark 2.17, since  $C$ is of full row rank,  $dim_{\mathbf{Q}}(R(C))=w.$ So $\sigma$ is an one-one linear mapping  from $\mathbf{Q}^w$ to $R(C)$.

The equation (36) is an identity in $x_1, x_2, \cdots,x_{w} .$
Each of the f's is a rational combination of the x's, since $f=x\, C$. By Remark 2.17, since  $C$ is of full row rank, each of the x's is a rational combination of the f's. Thus the equation (36) is an identity in
 the variables $f_1, f_2, \cdots,f_{w}, f_{w+1}, f_{w+2}$ for any  $f\in R(C)$.

 We express the integer $\alpha$ as the sum of four squares by Lemma 3.1, and bracket the terms $f_1^2+\cdots +f_{w}^2+f_{w+1}^2$  in fours.
 Each product of sums of four squares is itself a sum of four squares, and  so (36) yields
  $$\alpha(y_1^2+y_2^2+\cdots+y_{w}^2+y_{w+1}^2)+y_{w+2}^2  $$
\begin{equation}
=
\alpha(x_1^2+x_2^2+\cdots+x_{w}^2)+\beta z^2,
\end{equation}
 where $z=x_1+x_2+\cdots+x_{w}$, $f_{w+2}=y_{w+2}$, and the $y's$ are related to the
$f's$ by an invertible linear transformation with rational
coefficients.
 Thus $det(P)\neq 0$, $y=f \, P$.

Now define a linear mapping $\tau$ from $\mathbf{Q}^{w+2}$ to $\mathbf{Q}^{w+2}$
$$ \tau  : f\mapsto f\ P.$$
The image space $\tau \sigma(\mathbf{Q}^w)=V$ is a vector subspace of $\mathbf{Q}^{w+2}$ and $dim_{\mathbf{Q}}V=w$.

Since the $x's$ are rational linear combinations of
the $f's$, it follows that the $x's $ (and $z$) are rational linear
combinations of the $y's$.  Thus the equation (37) is an identity in
 the variables $y_1, y_2, \cdots,y_{w}, y_{w+1},y_{w+2}$  for any  $y\in V $.

Suppose that $ x_i=b_{i\, 1}y_1+\cdots +b_{i\, w}y_{w}+b_{i\, w+1}y_{w+1}+b_{i\, w+2}y_{w+2}$, $1\leq i\leq w$. We can
define $y_1$ as a rational linear combination of $y_2,\cdots ,
y_{w+1},y_{w+2}$, in such a way that $x_1^2=y_1^2$: if $b_{1\, 1}\neq 1$ we
set $y_1=\frac{1}{1-b_{1\, 1}}(b_{1\, 2}y_2+\cdots+b_{1\,
w+1}y_{w+1}+b_{1\, w+2}y_{w+2})$, while if  $b_{1\, 1}= 1$ we set
$y_1=\frac{1}{-1-b_{1\, 1}}(b_{1\, 2}y_2+\cdots+b_{1\,
w+1}y_{w+1}+b_{1\, w+2}y_{w+2})$. Now we know that  $x_2$ is  a rational linear
combination of the $y's$, and, using the relevant expression for
$y_1$ found above, we can express $x_2$  as a rational linear
combination of $y_2,\cdots,y_{w+1}, y_{w+2}$. As before, we fix $x_2$ as a
rational combination of $y_3,\cdots , y_{w+1}, y_{w+2}$ in such a way that
$x_2^2=y_2^2$. Continuing thus, we eventually obtain $x_1,\cdots
,x_{w}$ and $y_1,\cdots ,y_{w}$ as rational linear
combinations  of
$y_{w+1},y_{w+2}$, satisfying $x_i^2=y_i^2\, (1 \leq i \leq w)$.

We reduce the equations step by step in this way until a truncated triangle of equations is obtained, say
 $$y_1=d_{1\, 2}y_2+\cdots+d_{1\,
w+1}y_{w+1}+d_{1\,
w+2}y_{w+2},  $$
$$y_2=d_{2\, 3}y_3+\cdots+d_{2\,
w+1}y_{w+1}+d_{2\,
w+2}y_{w+2},  $$
$$\cdots \cdots\cdots\cdots\cdots   $$
$$y_w=d_{w\, \, w+1}y_{w+1}+d_{w\, \, w+2}y_{w+2};  $$
$$x_i^2=y_i^2, \, (1 \leq i \leq w); $$
where $d_{i\, j} \in \mathbf{ Q}$.

For any  $x\in \mathbf{Q}^w$ and $x\neq 0$,  by Remark 2.17, since  $C$ is of full row rank, $f=x\, C$, and $det(P)\neq 0$, $y=f \, P$, it implies $y\neq 0$. Let the last one $y_{w+2}\neq 0$ with suitable renumberings, if necessary.
Choose any non-zero rational value for  $y_{w+2} $.
All the
$x's $, the remaining $y's$, and $z$, are determined as above, and
substituting these values in (37) we obtain
\begin{equation}\alpha y_{w+1}^2+y_{w+2}^2= \beta z^2.
\end{equation}
Multiplying by a suitable constant  we have  that   $$ \alpha z^2= -x^2+ \beta y^2.$$  So the theorem is proved.  \pfend
\begin{exa}  The projective plane of order 7 is the symmetric $(57,8, 1)$ design.
Let $A$ be its incidence matrix, which is  a 57 by 57 matrix.
Choose its  bordered matrix $C$ is  a 59 by 61 matrix
as the following matrix.
 \[\begin{array}{ccccc} C=\left[\begin{array}{ccccc} A_{1\,\, 1} \ A_{1\,\, 2}
\\ A_{2\,\, 1}\ A_{2\,\, 2}
\end{array}\right ],
\end{array}\]
$$A_{1\,\, 1}=A .$$ $A_{1\,\, 2}$ is a 57 by 4 matrix  and
$$ A_{1\,\, 2}
=(1*1_{57}^t,  0*1_{57}^t,  0*1_{57}^t,  0*1_{57}^t). $$
 $A_{2\,\, 1}$ is
 a 2 by 57 matrix and
  \[\begin{array}{ccccc} A_{2\,\, 1}=\left[\begin{array}{ccccc}    0 \cdot1_{57} \\  0 \cdot1_{57}
\end{array}\right ].
\end{array}\]
$A_{2\,\, 2}$ is a 2 by 4 matrix and
 \[ \begin{array}{ccccc} A_{2\,\, 2}=\left[\begin{array}{ccccccccccccccc} 2\, \ \, 2  \  \,  1\,  \ 0 \\ 2 \ \, 0 \ -2  \   1 \end{array}\right ] . \end{array}
\]
  It is easy to check that $C$ has the property of the property of row inner products, i.e.,

  (i)\,\, the inner product of any two distinct rows of $C$ is equal to $2$ ;

  (ii)\,\, and
  the inner product of any rows with themselves of $C$ is equal to $9$.\linebreak
  It follows that
 $$ C\, C^t=7 I_{59}+2 J_{59}$$
  and $C$ is exactly  the  bordered matrix of the symmetric $(57,8, 1)$ design.
 We have that $\alpha=7, \beta=2$ and  the equation
  $$ 7 z^2=- x^2+ 2 y^2$$
must  have a solution in integers,
$x, y, z$,  not all zero  just as the assertion of the above theorem.
 \end{exa}

 \textbf{Proof of Main Theorem 1} \,\, By the above Lemmas we finish the proof of Main Theorem 1. \pfend
  Suppose there exists  a  $(v,k,\lambda)$ symmetric design with an incidence matrix $A$.
 It is difficult to construct a square  bordered  matrix of $A$. The author does this by computer computation in Maple.  But it is easy
 to construct a nonsquare  bordered matrix of $A$. The author also does this by computer computation in Maple just as the following remarks.

  \begin{rem}
   Suppose there exists  a  $(v,k,\lambda)$ symmetric design with an incidence matrix $A$.
  Further suppose there exists a positive integer $l$ such that $l=a^2+b^2$, where $a, b$ are two integers.
  By Lemma 3.7 and computation in Maple we can choose $l$ and construct  the bordered matrix $v+1$ by $v+2$ $C$ of  the
 incidence matrix $A$ as  the following matrix  \[\begin{array}{ccccc} C=\left[\begin{array}{ccccc} A_{1\,\, 1} \ A_{1\,\, 2}
\\ A_{2\,\, 1}\ A_{2\,\, 2}
\end{array}\right ],
\end{array}\]
$$A_{1\,\, 1}=A.$$
 $A_{1\,\, 2}$ is a $v$ by $2$ matrix  and $$ A_{1\,\, 2}
=( a\cdot 1_v^t,  b \cdot1_v^t).$$
 $A_{2\,\, 1}$ is
 a 1 by $v$ matrix and
  \[\begin{array}{ccccc} A_{2\,\, 1}=\left[\begin{array}{ccccc}  x_1 \cdot 1_{v}
\end{array}\right ],
\end{array}\]
 $A_{2\,\, 2}=(x_2,x_3)$ is some 1 by 2 matrix,  where $x_1, x_2,x_3$ are  some rational numbers.
\end{rem}
  \begin{rem}
   Suppose there exists  a  $(v,k,\lambda)$ symmetric design with an incidence matrix $A$.
  Further suppose there exists a positive integer $l$ such that $l=a^2+b^2$, where $a, b$ are two integers.
  By Lemma 3.7 and computation in Maple we can choose $l$ and construct  the bordered matrix $v+1$ by $v+3$ $C$ of  the
 incidence matrix $A$ as  the following matrix  \[\begin{array}{ccccc} C=\left[\begin{array}{ccccc} A_{1\,\, 1} \ A_{1\,\, 2}
\\ A_{2\,\, 1}\ A_{2\,\, 2}
\end{array}\right ],
\end{array}\]
$$A_{1\,\, 1}=A.$$
 $A_{1\,\, 2}$ is a $v$ by $3$ matrix  and $$ A_{1\,\, 2}
=( a\cdot 1_v^t,  b \cdot1_v^t, 0 \cdot1_v^t ).$$
 $A_{2\,\, 1}$ is
 a 1 by $v$ matrix and
  \[\begin{array}{ccccc} A_{2\,\, 1}=\left[\begin{array}{ccccc}  x_1 \cdot 1_{v}
\end{array}\right ],
\end{array}\]
 $A_{2\,\, 2}=(x_2,x_3,x_4)$ is some 1 by 3 matrix,  where $x_1, x_2,x_3, x_4$ are  some rational numbers.
\end{rem}

  \begin{rem}
   Suppose there exists  a  $(v,k,\lambda)$ symmetric design with an incidence matrix $A$.
  Further suppose there exists a positive integer $l$ such that $l=a^2+b^2$, where $a, b$ are two integers.
  By Lemma 3.7 and computation in Maple we can choose $l$ and construct  the bordered matrix $C$ of  the
 incidence matrix $A$ as  the following matrix  \[\begin{array}{ccccc} C=\left[\begin{array}{ccccc} A_{1\,\, 1} \ A_{1\,\, 2} \ A_{1\,\, 3}
\\ A_{2\,\, 1}\ A_{2\,\, 2} \ A_{2\,\, 3}
\end{array}\right ],
\end{array}\]
$$det(C\,\ C^t)\neq 0,$$
$$A_{1\,\, 1}=A.$$
 $A_{1\,\, 2}$ is a $v$ by $2$ matrix  and $$ A_{1\,\, 2}
=( a\cdot 1_v^t,  b \cdot1_v^t).$$
  $A_{1\, 3}$ is a $v$ by $s$ zero matrix, where $s$ is 1 or 2.   $A_{2\,\, 1}$ is
 a 2 by $v$ matrix and
  \[\begin{array}{ccccc} A_{2\,\, 1}=\left[\begin{array}{ccccc}  c \cdot 1_{v}
\\  d \cdot1_{v}
\end{array}\right ],
\end{array}\]
 where $c, d$ are some two rational numbers.
 $A_{2\,\, 2}$ is some 2 by 2 matrix.
$A_{2\,\, 3}$ is a 2 by $s$ matrix, where $s$ is 1 or 2.
\end{rem}
It is easy to construct  the above bordered matrix by the computer using Lemma 3.6, Lemma 3.7, Remark 3.8 and Remark 3.9   if it does exist just as In \S 5 and \S 6.

\begin{rem}
 Let $A$ be the incidence matrix of a  symmetric $(v,k,\lambda)$ design. Then the bordered  matrix of $A$ may not exist.
 If it exists   then  it is not unique for positive integers $s,l$.
 \end{rem}


\section{Proof of Main Theorem 2}

\begin{thm} Projective planes of order 6, 14, 21, 22, 30 and 33 do not exist.\end{thm}

\pf\,\ In these cases we can  use Theorem 2.11(the Bruck-Ryser Theorem).  \pfend
In order to use  main theorem 1 to show that symmetric designs with certain parameters cannot exist, we must show that the corresponding  bordered matrix exist and the corresponding equation has no integral solution.

It is easy to construct  the above bordered matrix by the computer using Lemma 3.6, Lemma 3.7, Remark 3.8 and Remark 3.9 in Maple   if it does exist. It should
be remarked that one does not need to trust the computer blindly. Although the proofs are discovered by the computer, it produces a proof certificate that can easily be checked by hand, if so desired.
\begin{thm}  There does not exist finite projective plane of order 10.
\end{thm}
\pf\,\, In this case we can not use the Bruck-Ryser Theorem but can use case 1 of Main Theorem 1.  Suppose that a  symmetric $(111,11, 1)$ design exists.
Let $A$ be its incidence matrix, which is  a 111 by 111 matrix.
Choose its  bordered matrix $C$ is  a 112 by 113 matrix
as the following matrix.
 \[\begin{array}{ccccc} C=\left[\begin{array}{ccccc} A_{1\,\, 1} \ A_{1\,\, 2}
\\ A_{2\,\, 1}\ A_{2\,\, 2}
\end{array}\right ],
\end{array}\]
$$A_{1\,\, 1}=A .$$ $A_{1\,\, 2}$ is a 111 by 2 matrix  and
$$ A_{1\,\, 2}
=(10*1_{111}^t,  0*1_{111}^t). $$
 $A_{2\,\, 1}$ is
 a 1 by 111 matrix and
  \[\begin{array}{ccccc} A_{2\,\, 1}=\left[\begin{array}{ccccc}    -\frac{2129}{11221}\cdot1_{111}
\end{array}\right ].
\end{array}\]
$A_{2\,\, 2}$ is a 1 by 2 matrix and
 \[ \begin{array}{ccccc} A_{2\,\, 2}=\left[\begin{array}{ccccccccccccccc}   \frac{115674}{11221} \ \frac{6}{7}   \end{array}\right ] . \end{array}
\]
It is easy to check that $C$ has the property of row inner products, i.e.,

  (i)\,\, the inner product of any two distinct rows of $C$ is equal to $101$ ;

  (ii)\,\, and
  the inner product of any rows with themselves of $C$ is equal to $111$.\linebreak
  It follows that
the property of  $$C\,\ C^t= 10 I_{112} + 101 J_{112}.$$ Thus
  $C$ is exactly  the  bordered matrix of the symmetric $(111,11, 1)$ design if $A$ exists. But by case 1 of Main Theorem  1 
  we have that 101  is  a  perfect square, 
 which is a contradiction. So there does not exist finite projective plane of order 10.
\pfend

\begin{thm}  There does not exist finite projective plane of order 12.
\end{thm}
\pf\,\, In this case we can not use the Bruck-Ryser Theorem but can use case 8 of Main Theorem 1. Suppose that a symmetric $(157,13, 1)$ design exists.
Let $A$ be its incidence matrix, which is  a 157 by 157 matrix.
Choose its  bordered  matrix $C$ is a 159 by 161 matrix
 as the following matrix.
 \[\begin{array}{ccccc} C=\left[\begin{array}{ccccc} A_{1\,\, 1} \ A_{1\,\, 2}  \ A_{1\,\, 3}
\\ A_{2\,\, 1}\ A_{2\,\, 2} \ A_{2\,\, 3}
\end{array}\right ],
\end{array}\]
$$A_{1\,\, 1}=A .$$
$A_{1\,\, 2}$ is a 157 by 2 matrix  and
$$ A_{1\,\, 2}
=(2*1_{157}^t, 0*1_{157}^t) .$$
$A_{1\,\, 3}$ is a 157 by 2 matrix  and
$$ A_{1\,\, 3}
=(0*1_{157}^t, 0*1_{157}^t) .$$
 $A_{2\,\, 1}$ is
 a 2 by 157 matrix and
  \[\begin{array}{ccccc} A_{2\,\, 1}=\left[\begin{array}{ccccc}  \frac{285}{2191 } \cdot 1_{157}
\\  -\frac{1}{7}\cdot1_{157}
\end{array}\right ].
\end{array}\]
 $A_{2\,\, 2}$ is a 2 by 2 matrix and
 \[ \begin{array}{ccccc} A_{2\,\, 2}=\left[\begin{array}{ccccccccccccccc}  \frac{3625}{2191} \ \frac{11}{7}    \\ \frac{24}{7} \ \frac{10}{7}   \end{array}\right ] . \end{array}
\]
$A_{2\,\, 3}$ is a 2 by 2 matrix and
 \[ \begin{array}{ccccc} A_{2\,\, 3}=\left[\begin{array}{ccccccccccccccc}  \frac{669}{313} \ \frac{669}{313}    \\ 0 \, \,  \  \, \,  0\, \,   \end{array}\right ] . \end{array}
\]
It is easy to check that $C$ has the property of row inner products, i.e.,

  (i)\,\, the inner product of any two distinct rows of $C$ is equal to $5$ ;

  (ii)\,\, and
  the inner product of any rows with themselves of $C$ is equal to $17$.\linebreak
  It follows that
the property of  $$C\,\ C^t= 12 I_{159} + 5 J_{159}.$$ Thus  $C$ is exactly  the  bordered matrix of the symmetric $(157,13, 1)$ design if $A$ exists. But by case 8 of  main Theorem 1  the equation
  $$ 12z^2=- x^2 +5  y^2$$
 must have a solution in integers,
$x, y, z$,  not all zero. It implies that, by Lemma 3.6, Lemma 3.7, Remark 3.8 and Remark 3.9,  the Legendre symbol $(\frac{5}{3})=1$, which is a contradiction. So there does not exist finite projective plane of order 12.

\pfend

\begin{thm}  There does not exist finite projective plane of order 15.
\end{thm}
\pf\,\, In this case we can not use the Bruck-Ryser Theorem but can use case 8 of  Main Theorem 1. Suppose that a symmetric $(241,16, 1)$ design exists.
Let $A$ be its incidence matrix, which is  a 241 by 241 matrix.
Choose its  bordered matrix $C$ is a 243 by 245 matrix
as the following matrix.
 \[\begin{array}{ccccc} C=\left[\begin{array}{ccccc} A_{1\,\, 1} \ A_{1\,\, 2} \ A_{1\,\, 3}
\\ A_{2\,\, 1}\ A_{2\,\, 2} \ A_{2\,\, 3}
\end{array}\right ],
\end{array}\]
$$A_{1\,\, 1}=A .$$
 $A_{1\,\, 2}$ is a 241 by 2 matrix  and
$$ A_{1\,\, 2}
=(7\cdot 1_{241}^t, 0*1_{241}^t ).$$
$A_{1\,\, 3}$ is a 241 by 2 matrix  and
$$ A_{1\,\, 3}
=(0*1_{241}^t, 0*1_{241}^t) .$$
$A_{2\,\, 1}$ is
 a 2 by 241 matrix and
  \[\begin{array}{ccccc} A_{2\,\, 1}=\left[\begin{array}{ccccc}  \frac{1432}{49911  } \cdot 1_{241}
\\  -\frac{23}{381}\cdot1_{241}
\end{array}\right ].
\end{array}\]
 $A_{2\,\, 2}$ is a 2 by 2 matrix and
 \[ \begin{array}{ccccc} A_{2\,\, 2}=\left[\begin{array}{ccccccccccccccc}  \frac{353234}{49911} \ -\frac{1}{3}    \\ \frac{2774}{381} \ \frac{10}{3}   \end{array}\right ] . \end{array}
\]
$A_{2\,\, 3}$ is a 2 by 2 matrix and
 \[ \begin{array}{ccccc} A_{2\,\, 3}=\left[\begin{array}{ccccccccccccccc}  \frac{120}{131} \ \frac{486}{131}    \\ 0 \, \,  \  \, \,  0\, \,   \end{array}\right ] . \end{array}
\]
It is easy to check that $C$ has the property of row inner products, i.e.,

  (i)\,\, the inner product of any two distinct rows of $C$ is equal to $50$ ;

  (ii)\,\, and
  the inner product of any rows with themselves of $C$ is equal to $65$.\linebreak
  It follows that the property of  $$C\,\ C^t= 15 I_{243} + 50 J_{243}.$$ Thus
  $C$ is exactly  the  bordered matrix of the symmetric $(241,16, 1)$ design if $A$ exists. But by case 8 of  Main Theorem 1  the equation
  $$ 15 z^2= -x^2 + 50 y^2$$
 must have a solution in integers,
$x, y, z$,  not all zero. It implies that, by Lemma 3.6, Lemma 3.7, Remark 3.8 and Remark 3.9,  the Legendre symbol $(\frac{2}{5})=1$, which is a contradiction. So there does not exist finite projective plane of order 15.

\begin{thm}  There does not exist finite projective plane of order 18.
\end{thm}
\pf\,\, In this case we can not use the Bruck-Ryser Theorem but can use case 1 of  Main Theorem 1. Suppose that a symmetric $(343,19, 1)$ design exists.
Let $A$ be its incidence matrix, which is a 343 by 343 matrix.
Choose its  bordered matrix $C$ is a 344 by 345 matrix
as the following matrix.
 \[\begin{array}{ccccc} C=\left[\begin{array}{ccccc} A_{1\,\, 1} \ A_{1\,\, 2}
\\ A_{2\,\, 1}\ A_{2\,\, 2}
\end{array}\right ],
\end{array}\]
$$A_{1\,\, 1}=A .$$
$A_{1\,\, 2}$ is a 343 by 2 matrix  and
$$ A_{1\,\, 2}
=(6\cdot 1_{343}^t, 0*1_{343}^t ).$$
$A_{2\,\, 1}$ is
 a 1 by 343 matrix and
  \[\begin{array}{ccccc} A_{2\,\, 1}=\left[\begin{array}{ccccc}   -\frac{23}{355}\cdot1_{343}
\end{array}\right ].
\end{array}\]
 $A_{2\,\, 2}$ is a 1 by 2 matrix and
 \[ \begin{array}{ccccc} A_{2\,\, 2}=\left[\begin{array}{ccccccccccccccc} 
 \frac{2262}{355} \ \frac{18}{5}   \end{array}\right ] . \end{array}
\]
It is easy to check that $C$ has  the property of row inner products, i.e.,

  (i)\,\, the inner product of any two distinct rows of $C$ is equal to $37$ ;

  (ii)\,\, and
  the inner product of any rows with themselves of $C$ is equal to $55$.\linebreak
  It follows that  the property of  $$C\,\ C^t= 18 I_{344} + 37 J_{344}.$$
Thus $C$ is exactly  the  bordered matrix of the symmetric $(343,19, 1)$ design if $A$ exists.
  But by case 1 of  Main Theorem 1 we have 37 is a perfect square,
 which is a contradiction. So there does not exist finite projective plane of order 18.
\pfend
\begin{thm}  There does not exist finite projective plane of order 20.
\end{thm}
\pf\,\,  In this case we can not use the Bruck-Ryser Theorem but can use case 8 of  Main Theorem 1. Suppose that a symmetric $(421,21, 1)$ design exists.
Let $A$ be its incidence matrix, which is  a 421 by 421 matrix.
Choose its  bordered matrix $C$ is  a 423 by 425 matrix
 as the following matrix.
 \[\begin{array}{ccccc} C=\left[\begin{array}{ccccc} A_{1\,\, 1} \ A_{1\,\, 2} \ A_{1\,\, 3}
\\ A_{2\,\, 1}\ A_{2\,\, 2} \ A_{2\,\, 3}
\end{array}\right ],
\end{array}\]
$$A_{1\,\, 1}=A .$$
$A_{1\,\, 2}$ is a 421 by 2 matrix  and
$$ A_{1\,\, 2}
=(8\cdot 1_{421}^t, 0*1_{421}^t ).$$
$A_{1\,\, 3}$ is a 421 by 2 matrix  and
$$ A_{1\,\, 3}
=(0*1_{421}^t, 0*1_{421}^t) .$$
$A_{2\,\, 1}$ is
 a 2 by 421 matrix and
  \[\begin{array}{ccccc} A_{2\,\, 1}=\left[\begin{array}{ccccc}  - \frac{549047}{30808125}   \cdot 1_{421}
\\  \frac{9231}{257419 }\cdot1_{421}
\end{array}\right ].
\end{array}\]
 $A_{2\,\, 2}$ is a 2 by 2 matrix and
 \[ \begin{array}{ccccc} A_{2\,\, 2}=\left[\begin{array}{ccccccccccccccc}  \frac{83919088}{10269375 } \  -\frac{2}{25}    \\ \frac{2067298}{257419} \ \frac{210}{47}   \end{array}\right ] . \end{array}
\]
$A_{2\,\, 3}$ is a 2 by 2 matrix and
 \[ \begin{array}{ccccc} A_{2\,\, 3}=\left[\begin{array}{ccccccccccccccc}  \frac{3808}{5625 } \ \frac{23614}{5625 }   \\ 0 \, \,  \  \, \,  0\, \,   \end{array}\right ] . \end{array}
\]
It is easy to check that $C$ has  the property of row inner products, i.e.,

  (i)\,\, the inner product of any two distinct rows of $C$ is equal to $65$ ;

  (ii)\,\, and
  the inner product of any rows with themselves of $C$ is equal to $85$.\linebreak
  It follows that the property of  $$C\,\ C^t= 20 I_{423} + 65 J_{423}.$$ Thus $C$ is exactly the  bordered
 matrix of the symmetric $(421,21, 1)$ design if $A$ exists. But by case 8 of  Main Theorem 1 the equation
  $$ 20z^2=-x^2+ 65y^2$$
 must have a solution in integers,
$x, y, z$,  not all zero.  It implies that , by Lemma 3.6, Lemma 3.7, Remark 3.8 and Remark 3.9,  the Legendre symbol $(\frac{13}{5})=1$, which is a contradiction. So there does not exist finite projective plane of order 20.

\begin{thm}  There does not exist finite projective plane of order 24.
\end{thm}
\pf\,\,  In this case we can not use the Bruck-Ryser Theorem but can use case 8 of  Main Theorem 1. Suppose that a symmetric $(601,25, 1)$ design exists.
Let $A$ be its incidence matrix, which is  a 601 by 601 matrix.
Choose its  bordered matrix $C$ is  a 603 by 605 matrix
 as the following matrix.
 \[\begin{array}{ccccc} C=\left[\begin{array}{ccccc} A_{1\,\, 1} \ A_{1\,\, 2}
\\ A_{2\,\, 1}\ A_{2\,\, 2}
\end{array}\right ],
\end{array}\]
$$A_{1\,\, 1}=A .$$
$A_{1\,\, 2}$ is a 601 by 4 matrix  and
$$ A_{1\,\, 2}
=(1\cdot 1_{601}^t, 0*1_{601}^t, 0*1_{601}^t, 0*1_{601}^t ).$$
$A_{2\,\, 1}$ is
 a 2 by 601 matrix and
  \[\begin{array}{ccccc} A_{2\,\, 1}=\left[\begin{array}{ccccc}   \frac{3}{25}   \cdot 1_{601}
\\  \frac{13}{185 }\cdot1_{601}
\end{array}\right ].
\end{array}\]
 $A_{2\,\, 2}$ is a 2 by 4 matrix and
 \[ \begin{array}{ccccc} A_{2\,\, 2}=\left[\begin{array}{ccccccccccccccc}  -1\ \frac{46}{25 } \  \frac{18}{5} \ 0    \\ \frac{9}{37} \ \frac{-284}{185} \ 0 \ \frac{168}{37}   \end{array}\right ] . \end{array}
\]
It is easy to check that $C$ has the property of row inner products, i.e.,

  (i)\,\, the inner product of any two distinct rows of $C$ is equal to $2$ ;

  (ii)\,\, and
  the inner product of any rows with themselves of $C$ is equal to $26$.\linebreak
  It follows that the property of  $$C\,\ C^t= 24 I_{603} + 2 J_{603}.$$
 Thus  $C$ is exactly the  bordered matrix of the symmetric $(601,24, 1)$ design if $A$ exists. But by case 8 of  Main Theorem 1 the equation
  $$ 24z^2=-x^2+ 2y^2$$
 must have a solution in integers,
$x, y, z$,  not all zero.  It implies that, by Lemma 3.6, Lemma 3.7, Remark 3.8 and Remark 3.9,  the Legendre symbol $(\frac{2}{3})=1$, which is a contradiction. So there does not exist finite projective plane of order 24.
\pfend
\begin{thm}  There does not exist finite projective plane of order 26.
\end{thm}
\pf\,\, In this case we can not use the Bruck-Ryser Theorem but can use case 1 of  Main Theorem 1. Suppose that a symmetric $(703,27, 1)$  design exists.
Let $A$ be its incidence matrix, which is a 703 by 703 matrix.
Choose its  bordered matrix $C$ is a 704 by 705 matrix
as the following matrix.
 \[\begin{array}{ccccc} C=\left[\begin{array}{ccccc} A_{1\,\, 1} \ A_{1\,\, 2}
\\ A_{2\,\, 1}\ A_{2\,\, 2}
\end{array}\right ],
\end{array}\]
$$A_{1\,\, 1}=A .$$
$A_{1\,\, 2}$ is a 703 by 2 matrix  and
$$ A_{1\,\, 2}
=(3\cdot 1_{703}^t, 0*1_{703}^t ).$$
$A_{2\,\, 1}$ is
 a 1 by 703 matrix and
  \[\begin{array}{ccccc} A_{2\,\, 1}=\left[\begin{array}{ccccc}   -\frac{4}{147}\cdot1_{703}
\end{array}\right ].
\end{array}\]
 $A_{2\,\, 2}$ is a 1 by 2 matrix and
 \[ \begin{array}{ccccc} A_{2\,\, 2}=\left[\begin{array}{ccccccccccccccc}
 \frac{526}{147} \ \frac{100}{21}   \end{array}\right ] . \end{array}
\]
  It is easy to check that $C$ has  the property of row inner products, i.e.,

  (i)\,\, the inner product of any two distinct rows of $C$ is equal to $10$ ;

  (ii)\,\, and
  the inner product of any rows with themselves of $C$ is equal to $36$.\linebreak
  It follows that the property of
  $$C\,\ C^t= 26 I_{704} + 10 J_{704},$$
   and $C$ is exactly  the  bordered matrix of the symmetric $(703,27, 1)$ design if $A$ exists.
  But by case 1 of  Main Theorem 1 we have that 10 is a perfect square,
 which is a contradiction. So there does not exist finite projective plane of order 26.
\pfend
\begin{thm}  There does not exist finite projective plane of order 28.
\end{thm}
\pf\,\,  In this case we can not use the Bruck-Ryser Theorem but can use case 8 of  Main Theorem 1. Suppose that a symmetric $(813,29, 1)$ design exists.
Let $A$ be its incidence matrix, which is  a 813 by 813 matrix.
Choose its  bordered matrix $C$ is  a 815 by 817 matrix
 as the following matrix.
 \[\begin{array}{ccccc} C=\left[\begin{array}{ccccc} A_{1\,\, 1} \ A_{1\,\, 2}
\\ A_{2\,\, 1}\ A_{2\,\, 2}
\end{array}\right ],
\end{array}\]
$$A_{1\,\, 1}=A .$$
$A_{1\,\, 2}$ is a 813 by 4 matrix  and
$$ A_{1\,\, 2}
=(1\cdot 1_{813}^t, 2*1_{813}^t, 0*1_{813}^t, 0*1_{813}^t ).$$
$A_{2\,\, 1}$ is
 a 2 by 813 matrix and
  \[\begin{array}{ccccc} A_{2\,\, 1}=\left[\begin{array}{ccccc}   \frac{1}{7}   \cdot 1_{813}
\\  \frac{291}{2590 }\cdot1_{813}
\end{array}\right ].
\end{array}\]
 $A_{2\,\, 2}$ is a 2 by 4 matrix and
 \[ \begin{array}{ccccc} A_{2\,\, 2}=\left[\begin{array}{ccccccccccccccc}   \frac{-23}{7 } \  \frac{18}{7} \ 0  \ 0  \\ \frac{5991}{2590} \ \frac{3}{14} \ \frac{336}{185} \ \frac{287}{74}   \end{array}\right ] . \end{array}
\]
 It is easy to check that $C$ has the property of row inner products, i.e.,

  (i)\,\, the inner product of any two distinct rows of $C$ is equal to $6$ ;

  (ii)\,\, and
  the inner product of any rows with themselves of $C$ is equal to $34$.\linebreak
  It follows that the property of  $$C\,\ C^t= 28 I_{815} + 6 J_{815}.$$
So $C$ is exactly the  bordered matrix of the symmetric $(813,29, 1)$ design if $A$ exists. But by case 8 of  Main Theorem 1 the equation
  $$ 28z^2=-x^2+ 6y^2$$
 must have a solution in integers,
$x, y, z$,  not all zero.  It implies that , by Lemma 3.6, Lemma 3.7, Remark 3.8 and Remark 3.9, the Legendre symbol $(\frac{6}{7})=1$, which is a contradiction. So there does not exist finite projective plane of order 28.
\pfend
\textbf{Proof of Main Theorem 2} \,\, By the above theorems we finish the proof of Main Theorem 2. \pfend

\section{Proof of  Main Theorem 3}
In order to use the main theorems to show that symmetric designs with certain parameters cannot exist, we must show that the corresponding  bordered matrix exist and the corresponding equation has no integral solution.
\begin{thm}  There does not exist symmetric $(49,16, 5)$ design.
\end{thm}
\pf\,\, In this case we can not use the Bruck-Ryser-Chowla Theorem but can use case 8 of  Main Theorem 1.  Suppose that a symmetric $(49,16, 5)$ design exists.
Let $A$ be its incidence matrix, which is  a 49 by 49 matrix.
Choose its  bordered matrix $C$ is  a 51 by 53 matrix
as the following matrix.
 \[\begin{array}{ccccc} C=\left[\begin{array}{ccccc} A_{1\,\, 1} \ A_{1\,\, 2}
\\ A_{2\,\, 1}\ A_{2\,\, 2}
\end{array}\right ],
\end{array}\]
$$A_{1\,\, 1}=A .$$ $A_{1\,\, 2}$ is a 49 by 4 matrix  and
$$ A_{1\,\, 2}
=(1*1_{49}^t,  0*1_{49}^t,0*1_{49}^t, 0*1_{49}^t ). $$
 $A_{2\,\, 1}$ is
 a 2 by 49 matrix and
  \[\begin{array}{ccccc} A_{2\,\, 1}=\left[\begin{array}{ccccc}    \frac{1}{3}\cdot 1_{49} \\  \frac{154}{425}\cdot 1_{49}
\end{array}\right ].
\end{array}\]
$A_{2\,\, 2}$ is a 2 by 4 matrix and
 \[ \begin{array}{ccccc} A_{2\,\, 2}=\left[\begin{array}{ccccccccccccccc}   \frac{2}{3}\,  \ \, \,  \frac{10}{3}\,  \  \, 0 \,  \  \,  0  \\   \frac{86}{425} \  \frac{-2}{125} \ \frac{242}{425} \ \frac{6787}{2125}  \end{array}\right ] . \end{array}
\]
  It is easy to check that $C$ has  the property of row inner products, i.e.,

  (i)\,\, the inner product of any two distinct rows of $C$ is equal to $6$ ;

  (ii)\,\, and
  the inner product of any rows with themselves of $C$ is equal to $17$.\linebreak
  It follows that  $$C\,\ C^t= 11 I_{51} + 6 J_{51},$$ and $C$ is exactly  the  bordered matrix of the symmetric $(49,16, 5)$ design if $A$ exists.
  But by case 8 of  main Theorem 1 the equation
  $$ 11z^2=-x^2+ 6 y^2$$
 must have a solution in integers,
$x, y, z$,  not all zero.  It implies that, by Lemma 3.6, Lemma 3.7, Remark 3.8 and Remark 3.9,  the Legendre symbol $(\frac{6}{11})=1$, which is a contradiction. So there does not exist symmetric $(49,16, 5)$ design.
\pfend
\begin{thm}  There does not exist symmetric $(154,18, 2)$ design.
\end{thm}
\pf\,\, In this case we can not use the Bruck-Ryser Theorem but can use case 8 of main Theorem 1.  Suppose that a symmetric $(154,18, 2)$ design exists.
Let $A$ be its incidence matrix, which is  a 154 by 154 matrix.
Choose its  bordered matrix $C$ is  a 155 by 157 matrix
as the following matrix.
 \[\begin{array}{ccccc} C=\left[\begin{array}{ccccc} A_{1\,\, 1} \ A_{1\,\, 2}
\\ A_{2\,\, 1}\ A_{2\,\, 2}
\end{array}\right ],
\end{array}\]
$$A_{1\,\, 1}=A .$$ $A_{1\,\, 2}$ is a 154 by 3 matrix  and
$$ A_{1\,\, 2}
=(1*1_{154}^t,  0*1_{154}^t,0*1_{154}^t ). $$
 $A_{2\,\, 1}$ is
 a 1 by 154 matrix and
  \[\begin{array}{ccccc} A_{2\,\, 1}=\left[\begin{array}{ccccc}    \frac{3}{20}\cdot1_{154}
\end{array}\right ].
\end{array}\]
$A_{2\,\, 2}$ is a 1 by 3 matrix and
 \[ \begin{array}{ccccc} A_{2\,\, 2}=\left[\begin{array}{ccccccccccccccc}   \frac{3}{10} \ \frac{47}{20} \ \frac{63}{20}   \end{array}\right ] . \end{array}
\]
  It is easy to check that $C$ has  the property of row inner products, i.e.,

  (i)\,\, the inner product of any two distinct rows of $C$ is equal to $3$ ;

  (ii)\,\, and
  the inner product of any rows with themselves of $C$ is equal to $19$.\linebreak
  It follows that $$C\,\ C^t= 16 I_{155} + 3 J_{155},$$  and $C$ is exactly  the  bordered matrix of the symmetric $(154,18, 2)$ design if $A$ exists.
  But by case 8 of  Main Theorem 1 the equation
  $$ 16z^2=-x^2+ 3y^2$$
 must have a solution in integers,
$x, y, z$,  not all zero.  It implies that, by Lemma 3.6, Lemma 3.7, Remark 3.8 and Remark 3.9,  the Legendre symbol $(\frac{-1}{3})=1$, which is a contradiction. So there does not exist symmetric $(154,18, 2)$ design.
\pfend

\begin{thm}  There does not exist symmetric  $(115,19, 3)$ design.
\end{thm}
\pf\,\, In this case we can not use the Bruck-Ryser Theorem but can use case 1 of  Main Theorem 1.  Suppose that a  symmetric $(115,19, 3)$ design exists.
Let $A$ be its incidence matrix, which is  a 115 by 115 matrix.
Choose its  bordered matrix $C$ is  a 116 by 117 matrix
as the following matrix.
 \[\begin{array}{ccccc} C=\left[\begin{array}{ccccc} A_{1\,\, 1} \ A_{1\,\, 2}
\\ A_{2\,\, 1}\ A_{2\,\, 2}
\end{array}\right ],
\end{array}\]
$$A_{1\,\, 1}=A .$$ $A_{1\,\, 2}$ is a 115 by 2 matrix  and
$$ A_{1\,\, 2}
=(3*1_{115}^t,  0*1_{115}^t ). $$
 $A_{2\,\, 1}$ is
 a 1 by 115 matrix and
  \[\begin{array}{ccccc} A_{2\,\, 1}=\left[\begin{array}{ccccc}    \frac{3}{7}\cdot1_{115}
\end{array}\right ].
\end{array}\]
$A_{2\,\, 2}$ is a 1 by 2 matrix and
 \[ \begin{array}{ccccc} A_{2\,\, 2}=\left[\begin{array}{ccccccccccccccc}   \frac{9}{7} \ \frac{16}{7}    \end{array}\right ] . \end{array}
\]
  It is easy to check that $C$ has the property of row inner products, i.e.,

  (i)\,\, the inner product of any two distinct rows of $C$ is equal to $12$ ;

  (ii)\,\, and
  the inner product of any rows with themselves of $C$ is equal to $28$.\linebreak
  It follows that  $$C\,\ C^t= 16 I_{116} + 12 J_{116},$$  and $C$ is exactly  the  bordered matrix of the  symmetric $(115,19, 3)$  design if $A$ exists.
  But by case 1 of  Main Theorem 1  we have that 12  is  a  perfect square,
 which is a contradiction.  So there does not exist any symmetric $(115,19, 3)$ design.
\pfend
\textbf{Proof of Main Theorem 3} \,\, By the above theorems we finish the proof of Main Theorem 3. \pfend

\section{Concluding remarks}
We conclude the discussion on block designs by mentioning the very short proof of the Bruck-Ryser-Chowla
theorem on the existence of symmetric block designs, which  is motivated at least in part by the matrix equation of set intersections\cite{ry3}.
Let $A$ be the incidence matrix of the symmetric $(v, k, \lambda)$ design.  Ryser  dealt only
with the case of symmetric $(v, k, \lambda)$ designs with $v$ odd. The criterion for $v$ even is
elementary.  He formed the following bordered matrix of order $v + 1$\cite{ry3}
   \[\begin{array}{ccccc} A^*=\left[\begin{array}{ccccc} A_{1\,\, 1} \ A_{1\,\, 2}
\\ A_{2\,\, 1}\ A_{2\,\, 2}
\end{array}\right ]
\end{array}\] \begin{equation}
\end{equation}
where $A_{1\,\, 1}=A $, $A_{1\,\, 2}$ is a column  vector $1_{v}^t$,  $A_{2\,\, 1}$ is a row vector $1_{v}$
 and $A_{2\,\, 2}=\frac{k}{\lambda}$. He also defined the following diagonal matrices $D$ and $E$ of order $v + 1$
$$ D = diag ( l,\cdots, 1,-\lambda),  E = diag(k-\lambda,\cdots, k-\lambda, -\frac{k-\lambda}{\lambda}). $$
Then it follows
  that the matrices $D, E$, and $A^*$ are
interrelated by the equation $$A^*DA^{*\,\ t}=E.$$
Thus the existence of the symmetric $ (v, k,\lambda)$ design implies that the diagonal matrices $D$
and $E$ of order $v + 1$ are congruent to one another over the field of rational
numbers. The remainder of the argument proceeds along standard lines and utilizes the Witt cancellation law. He just gave a new proof and did not obtain new necessary conditions on the existence of symmetric $(v,k, \lambda)$ designs.

In this paper we consider the  bordered matrix $C$ of the symmetric $(v, k, \lambda)$ design with preserving some row inner product property for some positive integer $l$, which is different from the above one, such that
 \[\begin{array}{ccccc} C=\left[\begin{array}{ccccc} A_{1\,\, 1} \ A_{1\,\, 2}
\\ A_{2\,\, 1}\ A_{2\,\, 2}
\end{array}\right ],
\end{array}\] \begin{equation}  C\, C^t=(k-\lambda) I_{w}+(\lambda +l) J_{w}.
\end{equation}
where $A_{1\,\, 1}=A $, $A_{1\,\, 2}$,   $A_{2\,\, 1}$ and $A_{2\,\, 2}$ are submatrices over $\mathbf{Q}$.

The matrix equation (40) is of fundamental importance. But it is difficult
to deal with this matrix equation in its full generality. In this paper $C$ maybe nonsquare matrix. The equation (40) implies
positive definite matrix  $(k-\lambda) I_{w}+(\lambda +l) J_{w}$ of order $w$ is quasi-congruent to the identity matrix of order $w+d$ with  plus $d$ over the field of rational
numbers. The equation (40) certainly contains much more information than (39). The difficulty lies in utilizing this information in an effective manner. So the  bordered matrix of $C$ of the symmetric $(v, k, \lambda)$ design, which  preserves  some row inner product property for some positive integer $l$, is just considered  more  property of $(0,1)$-matrix.  Let $d$ be the difference between the number of columns and the number of rows of $C$ in (40). If $d>2$, then we do not obtain the Diophantine equations of  Legendre type. Thus in this paper we just consider that $d$ is 1 or 2.
This has been the key breakthrough since 1950.

It was proved by a computer search that there does not exist any projective plane of order 10 by  Lam, C.W.H., Thiel, L. and Swiercz, S.
This is not the first time that a computer has played
 an important role in proving a theorem. A notable earlier example is the
 four-color theorem.
It is easy to construct  the above bordered matrix by the computer using Lemma 3.6, Lemma 3.7, Remark 3.8 and Remark 3.9 in Maple  in Theorem 5.1. It should
be remarked that one does not need to trust the computer blindly. Although the proofs are discovered by the computer, it produces a proof certificate that can easily be checked by hand, if so desired.  So we obtain a proof in the traditional mathematical
 sense for nonexistence of finite projective plane of order 10 and some other cases.

 For Problem 2.7 or Conjecture 2.12 Lam's algorithm\cite{l4} is an exponential time algorithm. But the proof of  main Theorem 1 is just the Ryser-Chowla
 elimination procedure  in \cite{ch}. Thus our algorithm is a polynomial time algorithm. It fully reflects  our algorithm high efficiency.


\end{document}